\documentclass[twoside, 11pt]{article}
\usepackage{indentfirst,bm,latexsym,amsfonts,amssymb,amsmath,mathrsfs}
\RequirePackage{ifthen,calc}
\usepackage{theorem}
\usepackage[dvips]{color}

 \catcode`@=11
 \def\@evenhead{\hbox to\textwidth{\footnotesize\rm\thepage \hfill
  {\it }}} 

 \def\@oddhead{\hbox to \textwidth{\footnotesize{\it
  On the barrier problem of branching random walk in a
time-inhomogeneous random environment  } \hfill\thepage}}

\catcode`@=12
\newcommand\ack{\section*{Acknowledgement}}

\newtheorem{thm}{\noindent Theorem}[section]
\newtheorem{lem}[thm]{\noindent Lemma}
\newtheorem{cor}[thm]{\noindent Corollary}
\newtheorem{prop}[thm]{\noindent Proposition}

\newtheorem{remark}[thm]{\noindent Remark}
\newtheorem{lemma}[thm]{Lemma}
\theoremheaderfont{\normalfont\bfseries} \theorembodyfont{\slshape}
{\theorembodyfont{\rmfamily}\newtheorem{ex}[thm]{\noindent Example}}
{\theorembodyfont{\rmfamily}}
{\theorembodyfont{\rmfamily}}
{\theorembodyfont{\rmfamily}}
{\theorembodyfont{\rmfamily}}
{\theorembodyfont{\rmfamily}}
{\theorembodyfont{\rmfamily}}
\usepackage[colorlinks,citecolor=blue,linkcolor=blue]{hyperref}

 \def\beqlb{\begin{eqnarray}}\def\eeqlb{\end{eqnarray}}
 \def\beqnn{\begin{eqnarray*}}\def\eeqnn{\end{eqnarray*}}
 \newcommand{\bgeqn}{\begin{equation}}
\newcommand{\edeqn}{\end{equation}}
\def\ra{\rightarrow}

 \numberwithin{equation}{section}
\def\qed{\hfill$\square$\smallskip}
\def\no{\nonumber}
\def\L{{\mathcal L}}

\def\ra{\rightarrow}
\def\iy{\infty}

\def\bfE{{\mathbb{E}}}

\def\mbfE{{\mathbf{E}}}
\def\mbr{{\mathbb{R}}}
\def\bfP{{\mathbb{P}}}

\def\mbfP{{\mathbf{P}}}
\def\bfR{{\mathbb{R}}}

\def\bfN{{\mathbb{N}}}
\def\1{{\mathbf{1}}}

\begin{document}
\title{\bf On the barrier problem of branching random walk in a time-inhomogeneous random environment
}
\author{You Lv\thanks{Email: lvyou@dhu.edu.cn },~~Wenming Hong\thanks{Corresponding author. Email: wmhong@bnu.edu.cn}
\\
\\ \small College of Science, Donghua University,
\\ \small Shanghai 201620, P. R. China.
\\ \small School of Mathematical Sciences $\&$ Laboratory of Mathematics and Complex Systems,
\\ \small Beijing Normal University, Beijing 100875, China
}
\date{}
\maketitle


\noindent\textbf{Abstract}: We introduce a random absorption barrier to a supercritical branching random walk with an i.i.d. random environment $\{\L_n\}$ indexed by time $n,$ i.e.,
in each generation, only the individuals born below the barrier can survive and reproduce.
The barrier is set as $\chi_n+an^{\alpha},$ where $a,\alpha$ are two constants and $\{\chi_n\}$ is a random walk determined by the random environment. We show that for almost surely $\L:=\{\L_n\},$
the time-inhomogeneous branching random walk with barrier will become extinct (resp., survive with positive probability) if $\alpha<\frac{1}{3}$ or $\alpha=\frac{1}{3}, a<a_c$ (resp., $\alpha>\frac{1}{3}, a>0$ or $\alpha=\frac{1}{3}, a>a_c$), where $a_c$ is a positive constant determined by the random environment. The rates of extinction when $\alpha<\frac{1}{3}, a\geq0$ and $\alpha=\frac{1}{3}, a\in(0,a_c)$ are also obtained. These extend the main results in A\"{\i}d\'{e}kon $\&$ Jaffuel (2011) and Jaffuel (2012) to the random environment case. The influence of the random environment has been specified.

\smallskip


\smallskip


\smallskip

\section{Introduction}
\subsection{Description of the model}

Branching random walk on $\mathbb{R}$ with an i.i.d. random environment in time is a natural extension of time-homogeneous branching random walk. It contains two levels of randomness. The first randomness comes from the random environment. 
Each realization of the random environment drives a time-inhomogeneous branching random walk, which is the second stage of randomness.
Compared with the time-homogeneous branching random walk, this model no longer requires the particles in different generations having the same reproduction law, instead admits the laws to vary from generation to generation according to the random environment.
Some relevant literature studying this model is listed in Section 1.3.

   We describe the model as follows. Let $(\Pi,\mathcal{F}_{\Pi})$ be a measurable space and $\Pi\subseteq\tilde{\Pi}:=\{\mathfrak{m}:\mathfrak{m}~\text{is~a~probability~measure~on~}V\},$ where $V:=\bfN\times\bfR\times\bfR\times\cdots.$ The random environment $\L$ is defined as an i.i.d. sequence of random
 variables $\{\L_1$,~$\L_2$,~$\cdots,\L_n,\cdots\}$, where $\L_1$ takes values in $(\Pi,\mathcal{F}_{\Pi})$.
 Let $\nu$ be the law of $\L$, then we call the product space $(\Pi^{\bfN}, \mathcal{F}_{\Pi}^{^{\bigotimes}\bfN}, \nu)$ the \emph{environment space}. For any realization $L:=\{L_1$,~$L_2$,~$\cdots,L_n,\cdots\}$ of $\L$, the time-inhomogeneous branching random walk driven by the environment $L$ is a process constructed as follows.

(1)~At time $0,$ an initial particle $\phi$ in generation $0$ is located at the origin.

(2)~At time $1,$ the particle $\phi$ dies and gives birth to $N(\phi)$ children who form the first generation. These children are located at $\zeta_i(\phi), 1\leq i\leq N(\phi),$ where the distribution of the random vector $X(\phi):=(N(\phi), \zeta_1(\phi),\zeta_2(\phi),\ldots)$ is $L_1.$ Note that the values $\zeta_i(\phi)$ for $i> N(\phi)$ do not play any role in our model. We introduce them only for convenience. For example, we can take $\zeta_i(\phi)=0$ for any $i> N(\phi).$

(3)~ Similarly, at generation $n+1,$ every particle $u$ alive at generation $n$ dies and gives birth to $N(u)$ children. If we denote $\zeta_i(u), 1\leq i\leq N(u)$ the displacement of the children with respect to their parent $u$, then $X(u):=(N(u), \zeta_1(u),\zeta_2(u),\cdots)$ is of distribution $L_{n+1}.$ We should emphasize that conditionally on any given environment $L,$ all particles in this system always behave independently. 

Conditionally on $\L,$ we write $(\Gamma,\mathcal{F}_{\Gamma}, \mbfP_{\L})$ for the probability space under which the
time-inhomogeneous branching random walk is defined. The probability $\mbfP_{\L}$ is conventionally called a {\it quenched law}.  We define the probability $\mathbf{P}:=\nu\bigotimes\mbfP_{\L}$ on the product space $(\Pi^{\bfN}\times\Gamma,\mathcal{F}_{\Pi}^{^{\bigotimes}\bfN}\bigotimes\mathcal{F}_{\Gamma})$ such that for any $F\in \mathcal{F}_{\Pi}^{^{\bigotimes}\bfN}, ~G\in\mathcal{F}_{\Gamma},$ we have
\begin{eqnarray}\label{APP}\mathbf{P}(F\times G)=\int_{\L\in F}\mbfP_{\L}(G)~d\nu(\L).\end{eqnarray}
The marginal distribution of probability $\mathbf{P}$ on $\Gamma$ is usually called an {\it annealed law}. The quenched law $\mbfP_{\L}$ can be viewed as the conditional probability of $\mathbf{P}$ given $\L.$ Throughout this paper, we consider the case $F=\Pi^{\bfN}.$ Hence without confusion we also denote the annealed law $\mathbf{P}$ and abbreviate $\mathbf{P}(\Pi^{\bfN}\times G)$ to~$\mathbf{P}(G).$ Moreover, we write $\mbfE_{\L}$ and~$\mathbf{E}$ for the corresponding expectation of $\mbfP_{\L}$ and~$\mathbf{P}$ respectively.

We denote by $\mathbf{T}$ the (random) genealogical tree of the process. For a given particle $u\in\mathbf{T}$ we write $V(u)\in\bfR$ for the position of $u$ and $|u|$ for the generation at which $u$ is alive. Then $(\mathbf{T}, V, \mbfP_{\L}, \mathbf{P})$ is called the {\it branching random walk in the time-inhomogeneous random environment $\L$} (BRWre). Especially, if there exists a $\iota\in\Pi$ such that $\mathbf{P}(\L_1=\iota)=1$ thus $\mathbf{P}(\L_i=\iota)=1, \forall i\in\bfN^+:=\{1,2,\cdots,n,\cdots\},$ which is conventionally called the {\it degenerate environment} (or {\it constant environment}), 
then the BRWre degenerates to the {\it time-homogeneous branching random walk} (BRW). Of course, one can describe the model by point process; see Mallein and Mi{\l}o\'{s} \cite{MM2016}.

\subsection{The barrier problem of BRW}

In this subsection, we will recall some progress for the barrier problem of BRW, i.e., the constant environment situation, and we use $\bfP$ and $\bfE$ to denote the probability and the corresponding expectation for the model without random environment (e.g., BRW, random walk) in the rest of the present paper.

In order to answer some questions about parallel simulations studied in Lubachevsky et al. \cite{LSW1989,LSW1990}, the barrier problem of BRW was first introduced in Biggins et al. \cite{BLSW1991}. The conclusion in \cite{BLSW1991} is closely related to the first order of the asymptotic behavior of
$$M_n:=\min\{V(u): u\in\mathbf{T}, |u|=n\},$$ i.e., the minimal displacement of the particles in the $n-$th generation.
Hammersley \cite{H1974}, Kingman \cite{K1975} and Biggins \cite{B1976} showed that (under some mild assumptions,) there is a finite constant $r$ such that \beqlb\label{HKB}\lim\limits_{n\rightarrow\infty}\frac{M_n}{n}=r,~~~ {\rm a.s.},\eeqlb which is the first order of BRW. Addario-Berry and Reed \cite{AR2009} and Hu and Shi \cite{HS2009} considered the second order of $M_n$. For example, \cite{HS2009} showed that 
\beqlb\label{Hushi}\lim\limits_{n\rightarrow\infty}\frac{M_n-rn}{\log n}=r_0, ~~{\rm in~ Probability},\eeqlb
where $r_{0}$ is a nonzero finite constant.
The weak convergence of $M_n-rn-r_0\log n$ can be found in A\"{\i}d\'{e}kon \cite{A2013}.

 We introduce some notations for a better understanding of the barrier problem. On the tree $\mathbf{T}$ we define a partial order $>$ such that $u>v$ if $v$ is an ancestor of $u$. We write $u\geq v$ if $u>v$ or $u=v$. We define an {\it infinite path} $u_{\infty}$ through $\mathbf{T}$ as a sequence of particles $u_{\infty}:=(u_i,i\in\bfN)$ such that $$\forall i\in\bfN,~~ |u_i|=i,~~ u_{i+1}>u_{i},~~ u_0=\phi~(\text{the initial particle}).$$ For any $i\leq |u|,$ we conventionally write $u_i$ for the ancestor of $u$ in generation $i.$ Let $\mathbf{T}_n:=\{u\in\mathbf{T}:|u|=n\}$ be the set of particles of generation $n$ and $\mathbf{T}_\infty$ the collection of all infinite paths through $\mathbf{T}.$

 The so-called ``barrier" is actually a function $\varphi:\bfN\ra\bfR.$ For any $u\in\mathbf{T},$ $u$ and all its descendants will be removed when $V(u)>\varphi(|u|).$ In other words, a particle in this system can survive only if all its ancestors and itself born below the barrier. For a BRW with a supercritical underlying branching process (i.e., $\bfE(\sum_{|u|=1}1)>1$),
 a natural question to consider when we add a barrier is whether the system still survives with positive probability. Define the event
 $$\mathcal{S}_0:=\{\exists u_{\infty}=(u_0,u_1,u_2, \ldots u_n, \ldots)\in \mathbf{T_{\infty}}, \forall i\in\bfN,  V(u_i)\leq \varphi(i)\},$$
 then $\bfP(\mathcal{S}_0)$ is the survival probability of BRW with barrier. In the light of \eqref{HKB} and \eqref{Hushi}, if
 the barrier function is set as $\varphi(i):=ri+ai^{\alpha},$ a series of predecessors' achievements is listed as follows.

Under some mild assumptions, Biggins et al. \cite{BLSW1991} showed that $\bfP(\mathcal{S}_0)>0$ when $\alpha=1, a>0$ and $\bfP(\mathcal{S}_0)=0$ when $\alpha=1, a<0.$

As a refined version of the above conclusion, Jaffuel \cite{BJ2012} showed that $\bfP(\mathcal{S}_0)>0$ when $\alpha=\frac{1}{3}, a>a_0$ and $\bfP(\mathcal{S}_0)=0$ when $\alpha=\frac{1}{3}, a<a_0,$ where $a_0$ is a positive constant. Obviously, this conclusion implies that $\bfP(\mathcal{S}_0)=0$ if $a=0.$

 Let $\mathcal{Y}_n:=\sharp\{u, |u|=n:\forall i\leq n, V(u_i)\leq \varphi(i)\}$ be the number of the survival populations in the $n$-th generation when we add the barrier $\varphi.$ Hence we have $\bfP(\mathcal{S}_0)=\lim\limits_{n\rightarrow\infty}\bfP(\mathcal{Y}_n>0).$ Note that $\bfP(\mathcal{S}_0)=0$ when $\alpha=1, a\leq 0.$ A\"{\i}d\'{e}kon and Jaffuel \cite{AJ2011} studied the extinction rate and showed that there are two finite negative constants $r_1,~r_2$ ($r_2$ depends on $a$) such that $\lim\limits_{n\rightarrow\infty}n^{-1/3}\log\bfP(\mathcal{Y}_n>0)=r_1$ when $a=0$ and $\lim\limits_{n\rightarrow\infty}n^{-1}\log\bfP(\mathcal{Y}_n>0)=r_2$ when $\alpha=1, a<0.$

For the case $\alpha=1, a>0,$ Gantert et al. \cite{GHS2011} gave the asymptotic behavior of $\bfP(\mathcal{S}_0)$ as $a\downarrow 0.$ They proved that $\lim\limits_{a\downarrow 0}\sqrt{a}\log\bfP(\mathcal{S}_0)=r_3,$ where the constant $r_3\in(-\infty, 0).$ Mallein \cite{M2017} obtained the same conclusion by an alternative proof.  Both \cite{GHS2011} and \cite{M2017} dealt with the problem (the asymptotic behavior of $\bfP(\mathcal{S}_0)$) in probabilistic approaches. Under a special case (assuming that
the branching law is binary branching and the random walk steps are bounded), this problem can be solved by an equation method (see \cite{BG2011}). 

The results mentioned above are all under the assumption that the associated random walk (derived from the celebrated many-to-one formula, see \cite[Theorem1.1]{S2015}) has finite variance. If the variance is infinite but the associated random walk is in the domain of
attraction of an $\alpha^*$-stable law, $\alpha^*\in(1,2)$, Liu and Zhang \cite{LZ2019} showed that there exists a constant $a^*_0$ depending on $\alpha^*$ such that $\bfP(\mathcal{S}_0)>0$ when $\alpha=\frac{1}{\alpha^*+1}, a>a^*_0$ and $\bfP(\mathcal{S}_0)=0$ when $\alpha=\frac{1}{\alpha^*+1}, a<a^*_0.$

We should explain that \cite{BLSW1991}, \cite{BJ2012}, \cite{AJ2011}, \cite{GHS2011} and \cite{LZ2019} all suppose that $r=0,$ which is an assumption of the so-called ``boundary case". Our statement above is essentially consistent with the original results in the boundary case according to the linear transformation in \cite{GHS2011}.

Kesten \cite{K1978}, Derrida and Simon \cite{DS2007,DS2008},  Harris and Harris \cite{HH2007} have studied the barrier problem of branching Brownian motion, which can be viewed as the continuous analog of BRW with barrier.

\subsection{The minimal displacement of BRWre}

Similar to the time-homogeneous case introduced in Section 1.2, the asymptotic behavior of the minimal displacement of BRWre is the theoretical basis for the barrier problem of BRWre. In this part, we list some conclusions about the minimal displacement of BRWre.
The model BRWre was first introduced in Biggins and Kyprianou \cite{BK2004}. Recall the definition of $M_n$ and the annealed law $\mbfP.$  Huang and Liu \cite{HL2014} proved that
there is a finite constant $d$ such that $\lim\limits_{n\rightarrow\infty}\frac{M_n}{n}=d,~ \mathbf{P}-$a.s.  \cite{HL2014} also obtained the large deviation principles for the counting measure about the population of the BRWre. Conclusions on the central limit theorem of the BRWre can be found in Gao, Liu, Wang \cite{GLW2014} and Gao, Liu \cite{GL2016}. The moderate deviation principles and the $L^p$ convergence rate have been investigated in Wang and Huang \cite{WH2017}.

What inspires our work most is the second order of the asymptotic behavior of $M_n$ condiered in Mallein and Mi{\l}o\'{s} \cite{MM2016}. They showed that there exists a random walk $\{\chi_n\}_{n\in\bfN}$ (the precise expression of $\chi_n$ is given in \eqref{lmost}) with i.i.d. increments under the annealed law $\mbfP$ such that \beqlb\label{SOBRWre}\frac{M_n-\chi_n}{\log n}\ra c, ~~~n\ra\iy,~~~\text{in~~Probability}~~ \mathbf{P},\eeqlb
where $c$ is a finite constant and $\mathbf{E}(\chi_1)=d.$ 
$\eqref{SOBRWre}$ shows that for BRWre, the trajectory of
$\{M_n\}_{n\in\bfN}$ is around the random walk $\{\chi_n\}_{n\in\bfN}$ (but not $\{dn\}_{n\in\bfN}$) with a logarithmic correction, which is different from the corresponding behavior of BRW (see \eqref{Hushi}). This fact provides a helpful guidance on how to set a reasonable barrier in the random environment case.

Some other types of inhomogeneous branching random walk have been studied. For example, Mallein \cite{M2015b} studied the maximal displacement of a branching random walk in time-inhomogeneous but non-random environment. Baillon, Cl\'{e}ment, Greven and den Hollander \cite{BC1993} considered a branching random walk with the environment determined by the space instead of the time. Hu and Yoshida \cite{HY2009}, as well as other authors, took interest in the branching random walk in a space-time random environment.

\subsection{Structure of this paper}

The rest of the paper is organized as follows.

In Section 2, we introduce some notations, assumptions and the main results. Moreover, we give an example satisfying our assumptions.

The many-to-one formula of time-inhomogeneous bivariate version is given in Section 3.

Section 4 is devoted to the small deviation principle of a random walk with time-inhomogeneous random environment. 
This principle is a basic tool to solve the barrier problem of BRWre. 

In Section 5, we give the proof of the propositions and example stated in Section 2.

At last, we prove the main theorems in Section 6 by all the preparations in Sections 3-5.

We comment that in the present paper we not only give the new results (in Section 2) and the proof (in Sections 3-6), but also give a detail analysis and comprehensive interpretation in Sections 2 (the proofs to support our analysis are given in Sections 4 and 5) on the setting of the assumptions in our main results, including the origin, necessity, alternatives, substitutes, comparison, particularity of these assumptions and an example to satisfy the assumptions.
Since the assumptions look different and more complicated than the corresponding assumptions for the barrier problem of BRW, we take several pages to show that our assumptions are totally acceptable and reasonable.

\section{Basic assumptions and Main results}
\subsection{Notations and assumptions}

First, we give some notations for the model BRWre. For every $n\in\bfN^+,$ define the log-Laplace transform function
$$\kappa_n(\theta):=\log\mbfE_{\L}\left(\sum^{N(u)}_{i=1}e^{-\theta \zeta_i(u)}\right),~|u|=n-1,~ \theta\in[0,+\infty).$$
Note that it is well-defined since conditionally on any realization $L_n$ of $\L_n,$ for each $u$ in generation $n-1,$ $X(u)$ has the common law $L_n$.
Meanwhile we should note that for a fixed $\theta$, $\kappa_n(\theta)$ is a random variable determined by the random environment $\L$. (More precisely, it is determined by $\L_n$.) Therefore, $\{\kappa_n(\theta), n\in\bfN^+\}$ is a sequence of i.i.d. random variables since $\{\L_n, n\in\bfN^+\}$ is i.i.d.

Throughout the present paper, we assume that $\mathbf{E}(\kappa^-_n(\theta))<+\infty$ for all $\theta\geq 0,$ where $\kappa^-_n(\theta)$ is defined as $\max\{0,-\kappa_n(\theta)\}.$  Hence we can define $\kappa: [0, +\infty)\rightarrow (-\infty,+\infty]$ by $$\kappa(\theta):=\mbfE(\kappa_n(\theta)).$$

Now we introduce four basic assumptions in the present paper.
The time-homogeneous versions of Conditions 1 and 3 are classic assumptions which often appear in the papers on BRW (for example, the assumptions (1.3), (2.2)-(2.4) in \cite{GHS2011}).
Conditions 2 and 4 are automatically satisfied under Conditions 1 and 3 when the random environment degenerates.

{\rm Condition 1} {\it Assume that $\kappa(0)>0,$ there exists $0<\vartheta<\bar{\theta}$ such that $\kappa(\bar{\theta})<+\infty$ and  
\beqlb\label{sect-2}\kappa(\vartheta)=\vartheta\kappa'(\vartheta).\eeqlb}

{\rm Condition 2} {\it There exist constants $\lambda_1>3, \lambda_2>2$ such that
\begin{eqnarray}\label{c2'a}
\mathbf{E}\left(|\kappa_1(\vartheta)-\vartheta \kappa'_1(\vartheta)|^{2\lambda_1}\right)<+\infty;
\end{eqnarray}
\begin{eqnarray}\label{c2'b}
\mbfE\left(\left[\frac{\mbfE_{\L}\left(\sum_{i=1}^{N(\phi)}|\zeta_i(\phi)+\kappa'_1(\vartheta)|^{\lambda_2}e^{-\vartheta \zeta_i(\phi)}\right)}{\mbfE_{\L}\left(\sum_{i=1}^{N(\phi)}e^{-\vartheta \zeta_i(\phi)}\right)}\right]^{\lambda_1}\right)<+\infty.
\end{eqnarray}
}

{\rm Condition 3} {\it Assume that we can find constants $\lambda_3>6, \lambda_4>0$ such that \begin{eqnarray}\label{c3a}\mbfE(|\kappa_1(\vartheta+\lambda_4)|^{\lambda_3})+\mbfE(|\kappa_1(\vartheta)|^{\lambda_3})<+\infty,~~ \mbfE([\log^+\mbfE_{\L}(N(\phi) ^{1+\lambda_4})]^{\lambda_3})<+\infty,\end{eqnarray}
where $\log^+\mbfE_{\L}(N(\phi)^{1+\lambda_4})=\log\max\{\mbfE_{\L}(N(\phi)^{1+\lambda_4}),1\}.$ }

{\rm Condition 4} {\it There exists constants $\lambda_5> 2, y<0$ such that
\begin{eqnarray}\label{c5}\mathbf{E}\left(\left|\log\mbfE_{\L}\left(\sum_{i=1}^{N(\phi)}\1_{\{\vartheta \zeta_i(\phi)+ \kappa_1(\vartheta)\leq y\}} \right)\right|^{\lambda_5}\right)<+\infty.\end{eqnarray}}

At least for the method we use, Conditions 1-4 are the almost tight version (i.e. it can not be further weakened). Condition 1 is set to ensure that $M_n/n$ has a finite limit (all investigations about the barrier problems are based on this behavior of $M_n/n$ in the time-homogeneous case) and the underlying branching process in random environment is supercritical. A standard strategy to research
 the barrier problems in the time-homogeneous case is combining the small deviation principle of random walk with the second moment method. In the present paper, the outline of the strategy is as follows. Condition 3 is set for applying the second moment method and Conditions 2 and 4 for applying the small deviation principle given in Section 4, where Condition 4 is set specially for applying the lower bound of the small deviation principle.
Some propositions will be given for a better understanding of Conditions 1-4.
\subsubsection{Some explanations of Conditions 1-4}

In this section, for a better and intuitive understanding of the above conditions, we give some properties and descriptions on the conditions.
\begin{prop}\label{propcon1}
Condition 1 implies that $\mathbf{P}(\kappa''_1(\vartheta)>0)>0,$ hence $\kappa''(\vartheta)>0.$
\end{prop}
~~For the time-homogeneous case, this proposition is obvious, while for the random environment case it needs to be proved (see Section 5). Moreover, this proposition is a indispensable preparation for the proof of the main results (Theorems \ref{1} and \ref{2}). The following two propositions give some sufficient conditions (which has a more intuitive description) for Conditions 2-3.
\begin{prop}\label{propcon2}
If there exists $\lambda_6>\frac{3}{2}$ such that
\begin{eqnarray}\label{c21a}\exists\lambda_6>\frac{3}{2}, 
 ~\mbfE((\kappa^{(4)}_1(\vartheta)+3[\kappa''_1(\vartheta)]^2)^{\lambda_6})<+\infty, \end{eqnarray}  where $\kappa^{(n)}_1(\vartheta):=\frac{d^n\kappa(\theta)}{d\theta}|_{\theta=\vartheta},$ then $\eqref{c2'b}$ holds.

Furthermore, if there exist $\lambda_7,\lambda_8>0$ such that
\begin{eqnarray}\label{c22a}  \mbfE(e^{\kappa_1(\vartheta-\lambda_7)-\kappa_1(\vartheta)}+e^{\kappa_1(\vartheta+\lambda_8)-\kappa_1(\vartheta)})<+\infty,\end{eqnarray}
then \eqref{c21a} holds.
\end{prop}
\begin{prop}\label{propcon3}
If there exist $\lambda_9,\lambda_{10},\lambda_{11}>0,\lambda_{12}>6$ such that
\begin{eqnarray}\label{p3}\mbfE(e^{\kappa_1(\vartheta-\lambda_{9})-\kappa_1(\vartheta)}+e^{|\kappa_1(\vartheta+\lambda_{10})-\kappa_1(\vartheta)|}+|\kappa_1(\vartheta)|^{\lambda_{12}})+\mbfE (N(\phi)^{1+\lambda_{11}})<+\infty,\end{eqnarray}
then Conditions 2-3 hold.
\end{prop}

Though \eqref{p3} looks like more intuitive than Condition 2 and Condition 3, we should also note that it is more restrictive than Conditions 2-3 since in \eqref{c2'a}, \eqref{c2'b} and \eqref{c3a}, we actually do not need $\kappa_1(\vartheta-\lambda_{9})-\kappa_1(\vartheta), |\kappa_1(\vartheta+\lambda_{10})-\kappa_1(\vartheta)|$ and $\log^{+}\mbfE_{\L}(N(\phi)^{1+\lambda_{11}})$ have finite exponential moments. We remind that Mallein and Mi{\l}o\'{s} \cite{MM2016} (which considered the minimal displacement of BRWre, see \eqref{SOBRWre}) requires that $\kappa_1(\vartheta-\lambda_{9})-\kappa_1(\vartheta), |\kappa_1(\vartheta+\lambda_{10})-\kappa_1(\vartheta)|, \log^{+}\mbfE_{\L}(N(\phi)^{1+\lambda_{11}})$ have finite exponential moments (see \cite[(1.7)-(1.8)]{MM2016}). More exactly, the main results in the present paper do not need the associated random walk $\{T_n\}_{n\in\bfN}$ defined in \eqref{mtoT} has finite exponential moment while \eqref{SOBRWre} needs that.

An explanation of Condition 4 is given in Remark \ref{5}, which shows that it is a necessary condition and specially needed when we consider the barrier problem in some extent. (We point out that Mallein and Mi{\l}o\'{s} \cite{MM2016} do not need this assumption since it does not involve the barrier problem.)

A comparison between Conditions 1-4 and the corresponding assumptions for the barrier problem of BRW is given in Remark \ref{compare}, in which we can trace why we set these conditions in this way and see some difficulties brought from the random environment.

Since the proofs of Propositions \ref{propcon2}-\ref{propcon3} will involve the many-to-one formula and the associated random walk, which are introduced in Sections 3 and 4 respectively, we might as well put all the proofs of the propositions after Section 4.

\subsubsection{An example which satisfies Conditions 1-4}

Recalling the definition of BRWre in Section 1 we see that the law of $\L_1$ totally determines a BRWre.
~Denote $X(\mathfrak{m}):=(N(\mathfrak{m}),\zeta_1(\mathfrak{m}),\zeta_2(\mathfrak{m}),\ldots)$ a random vector with distribution $\mathfrak{m}\in\Pi$ ($\Pi$ have been defined in Section 1.1) taking values on $\bfN\times\bfR\times\bfR\times\bfR\cdots.$
 Next we give an example in which the branching law and the displacement law are independent of each other.
\begin{ex}\label{exb0}
Let $\mu(\cdot)$ and $\sigma(\cdot)$ be two functionals defined on $\Pi$ and $\sigma(\mathfrak{m})>0, \forall \mathfrak{m}\in\Pi$. For any $\mathfrak{m}\in\Pi,$ we assume that for any realization of $N(\mathfrak{m})$, $(\zeta_i(\mathfrak{m}), i\leq N(\mathfrak{m}))$ is a sequence of i.i.d. random variables with a common normal distribution $\mathcal{N}(\mu(\mathfrak{m}),\sigma^2(\mathfrak{m})),$ where $\sigma^2(\mathfrak{m}):=(\sigma(\mathfrak{m}))^2.$
Here we remind that the values of $\mu(\mathfrak{m})$ and $\sigma(\mathfrak{m})$ are not affected by the realization of $X(\mathfrak{m}).$
We assume that there exist two constants $\tau_1>6, \tau_2>4$ such that
\begin{eqnarray}\label{exb1}\mbfE(\log\mbfE_{\L}N(\mathfrak{m}))>0, ~~\mathbf{E}[(\log^+N(\mathfrak{m}))^{\tau_1}]<+\infty,\end{eqnarray}
\begin{eqnarray}\label{exb2}\mathbf{E}(\sigma(\mathfrak{m})^{2\tau_1})<+\infty,~~ \mathbf{E}(\sigma(\mathfrak{m})^{-\tau_2})<+\infty.\end{eqnarray}
Then this example satisfies Conditions 1-4.
\end{ex}

Note that in this example we have not any requirement on the random variable $\mu(\mathfrak{m}).$

 The proof of this example also involves the many-to-one formula and the associated random walk hence we also put the proof after Section 4.



\subsection{Main results}

First we introduce the barrier considered in the present paper. The enlightenment about how to set the barrier function is from the main result in Mallein and Mi{\l}o\'{s} \cite{MM2016}.
Recall that $V(u)$ presents the position of particle $u$. Mallein and Mi{\l}o\'{s} \cite{MM2016} have shown that \beqlb\label{lmost}\lim\limits_{n\rightarrow+\infty}\frac{\min_{|u|=n}V(u)+\vartheta^{-1}K_n}{\log n}=c,~~~~~\text{in Probability}~\mbfP,\eeqlb
where $c$ is a finite constant and \beqlb\label{Kn}K_n:=\sum_{i=1}^n\kappa_i(\vartheta),~~~~ K_0:=0.\eeqlb
Hence we see for BRWre, the time-homogeneous random walk $\{\vartheta^{-1}K_n\}_{n\in\bfN}$ gives the first order of the asymptotic behavior of the minimal displacement. (Of course, this random walk is exactly the $\{\chi_n\}_{n\in\bfN}$ mentioned in \eqref{SOBRWre}). Recall that the first order of the minimal displacement of a BRW is $rn$ ($r$ is the one introduced in Section 1.2) and hence the barrier function in \cite{AJ2011}, \cite{GHS2011} and \cite{BJ2012} is set as the form $\varphi(i):=ri+ai^{\alpha}.$ Through the above analysis, in the present paper we set the
barrier function $\varphi_{\L}(i):=-\vartheta^{-1}K_i+ai^{\alpha},$ where we use the notation $\varphi_{\L}(i)$ but not $\varphi(i)$ since $\L$ brings randomness to $K_n.$ Hence we can see that the setting of the barrier is an important difference between our main results and the achieved results on the barrier problem of BRW. In the present paper, we set the random barrier $\varphi_{\L}(i)$ according to the random environment while for BRW, the barrier $\varphi(i)$ is non-random because of the constant environment. Moreover, the barrier $\varphi_{\L}(i)$ we set will become to $\varphi(i)$ when the random environment is degenerate.

Before giving the main results, we need to introduce the function $\gamma$ with the definition 
\beqlb\label{bmt}\gamma(\beta):=\lim\limits_{t\rightarrow+\infty}\frac{-\log \mbfP(\forall_{s\leq t} B_s\in[-\frac{1}{2}+\beta W_s, \frac{1}{2}+\beta W_s]|W)}{t},~~~\rm{a.s.,}\eeqlb where $B, W$ are two independent standard Brownian motions and $B_0=0, W_0=0$. $\gamma$ has been introduced in Lv \cite[Theorem 2.1]{LY201801}. \cite{LY201801} showed that for any given realization of $W$ (in the sense of almost surely), $\gamma$ is a well-defined, positive, convex and even function with $\gamma(0)=\frac{\pi^2}{2}.$
 Moreover, $\gamma$ is strictly increasing on $[0,+\infty)$ and strictly decreasing on $(-\infty,0].$

 Throughout this paper, we denote
\beqlb\label{AQ}\sigma_A:=\sqrt{\mathbf{E}\left(\Big(\kappa_1(\vartheta)-\vartheta\kappa'_1(\vartheta)\Big)^2\right)},~~\sigma_Q:=\vartheta\sqrt{\mathbf{E}(\kappa''_1(\vartheta))},~~\gamma_\sigma:=\sigma^2_Q\gamma\left(\frac{\sigma_A}{\sigma_Q}\right).\eeqlb
The following two theorems are the main results in the present paper.
\begin{thm}\label{1} $(${\bf critical criterion}$)$~Let the barrier function be $\varphi_{\L}(i):=-\vartheta^{-1}K_i+a i^{\alpha},~\forall i\in\bfN.$
Recall the definition of the infinite path $u_{\infty}$ and their collection $\mathbf{T_{\infty}}.$ Define the event $$\mathcal{S}:=\{\exists u_{\infty}:=(u_1,u_2, \ldots u_n, \ldots)\in \mathbf{T_{\infty}}, \forall i\in\bfN,  V(u_i)\leq \varphi_{_\L}(i)\},$$
which is the event that the system still survives after we add the barrier $\varphi_{_\L}.$
Denote $a_c:=\frac{3\sqrt[3]{6\gamma_\sigma}}{2\vartheta}.$
\begin{itemize}
\item Suppose that Conditions 1-4 hold, 
then the following two statements are true.

{\rm (1a)}~If~$\alpha>\frac{1}{3}, a>0$, then~$\mbfP_{\L}(\mathcal{S})>0,~{\rm \mathbf{P}-a.s.}$

{\rm (1b)}~If~$\alpha=\frac{1}{3}, a>a_c$, then~$\mbfP_{\L}(\mathcal{S})>0,~{\rm \mathbf{P}-a.s.}$

\item Suppose that Conditions 1-2 hold, then the following two statements are true.

{\rm (1c)}~If~$\alpha=\frac{1}{3}, a<a_c$, then~$\mbfP_{\L}(\mathcal{S})=0,~{\rm \mathbf{P}-a.s.}$

{\rm (1d)}~If~$\alpha<\frac{1}{3}, a\in\bfR$, then~$\mbfP_{\L}(\mathcal{S})=0,~{\rm \mathbf{P}-a.s.}$
\end{itemize}
\end{thm}

\begin{thm}\label{2} $(${\bf extinction rate}$)$ Define $$Y_n:=\sharp\{|u|=n:~\forall i\leq n, ~V(u_i)\leq \varphi_{\L}(i)\},$$
which represents the number of the surviving particles in generation $n$ after we add the barrier $\varphi_{\L}$.
\begin{itemize}\item
Suppose that Conditions 1-3 hold and \eqref{c5} holds with some $\lambda_5\geq1$, then we have the following {\rm (2a)-(2c)}.

{\rm (2a)}. Let $\alpha=\frac{1}{3}, a\in(0,a_c).$ Then
\begin{eqnarray}\label{2a}
\lim\limits_{n\rightarrow\infty}\frac{\log\mbfP_{\L}(Y_n>0)}{\sqrt[3]{n}}=-\vartheta\tilde{q}(0),~{\rm \mathbf{P}-a.s.},
\end{eqnarray}
where $\tilde{q}$ is the unique solution in $\mathcal{C}[0,1]$ of the integral equation
\begin{eqnarray}\label{2a1}\forall ~t\in[0,1],~~ -\vartheta{q}(0)=\vartheta at^{\frac{1}{3}}
 -\vartheta {q}(t)-\gamma_\sigma\int_0^{t}(\vartheta {q}(x))^{-2}dx \end{eqnarray} and satisfies
\begin{eqnarray}\label{2a2}\tilde{q}(0)>0,~\tilde{q}(1)=0,~\int_0^1\frac{1}{(\tilde{q} (x))^2}dx<+\infty.\end{eqnarray}

{\rm (2b)}. Let $\alpha\in(0,\frac{1}{3}), a\geq0.$ 
Then it is true that\begin{eqnarray}\label{2b}
\lim\limits_{n\rightarrow\infty}\frac{\log\mbfP_{\L}(Y_n>0)}{\sqrt[3]{n}}=-\sqrt[3]{3\gamma_{\sigma}},~{\rm \mathbf{P}-a.s.}
\end{eqnarray}

{\rm (2c)}. Define
$$P_{\L}(n,c):=\mbfP_\L\left(\exists |u|= n, \forall i\leq n, V(u_i)\leq ci-\vartheta^{-1}K_{i}\right).$$
For any constant $b>0,$ we have
\begin{eqnarray}\label{2c1} \varlimsup\limits_{n\rightarrow\infty}\frac{1}{n^{1/3}}\log P_{\L}\left(n,bn^{-\frac{2}{3}}\right)\leq -x_b,~~{\rm~~\mathbf{P} -a.s.},\end{eqnarray}
\begin{eqnarray}\label{2c2}\varliminf\limits_{n\rightarrow\infty}\frac{1}{n^{1/3}}\log P_{\L}\left(n,bn^{-\frac{2}{3}}\right)\geq -\sqrt{\frac{\gamma_{\sigma}}{\vartheta b}},~~{\rm~~\mathbf{P} -a.s.},\end{eqnarray}
where $x_b$ is the solution of ~$\frac{3\gamma_{\sigma}}{x^2}-x=3\vartheta b$ on $(0,+\infty).$
\item Suppose that Conditions 1-4 hold, then we have the following {\rm (2d)}.

{\rm (2d)}. Let~$\alpha=\frac{1}{3},~ a>a_c.$ For any given constant $\varepsilon>0,$ there exists $M\in\bfN$ large enough such that
\begin{eqnarray}\label{big1}
\mbfP_{\L}\left(\varliminf\limits_{k\rightarrow+\infty}\frac{\log Y_{M^k}}{M^{\frac{k}{3}}}>b_2\vartheta-\varepsilon\right)>0,~{\rm \mathbf{P}-a.s.},\end{eqnarray}
where $b_2$ is the maximum $b$ satisfying~$\vartheta a=\vartheta b+\frac{3\gamma_{\sigma}}{b^2\vartheta^2}$.
\end{itemize}
\end{thm}\begin{remark}
 If the random environment degenerates to time-homogeneous environment (i.e., the model BRWre degenerates to BRW), then Condition 1 implies that $\sigma_A=0,$ and hence $\gamma_{\sigma}=\sigma^2_Q\gamma(0)=\frac{\pi^2\sigma^2_Q}{2}.$
~That is to say, our results are consistent with the corresponding conclusions about BRW in \cite{AJ2011} and \cite{BJ2012}.
\end{remark}

\begin{remark}\label{compare}
Let us compare Conditions 1-4 and the corresponding assumptions for the barrier problem of BRW in \cite{BJ2012}. When the model BRWre degenerates to BRW, Condition 1 and Condition 3 are
the basic assumptions in \cite{BJ2012}, which considers the barrier problem of BRW. In our notations, the basic assumptions in \cite{BJ2012} can be stated as follows: there exist constants $\vartheta, \lambda_{13}, \lambda_{14}>0$ such that
\begin{eqnarray}\label{Jafassume}\kappa_1(\vartheta)-\vartheta\kappa'_1(\vartheta)=0,~\Phi_1(\vartheta+\lambda_{13})<+\infty,~ \mbfE_{\L}(N(\phi)^{1+\lambda_{14}})<+\infty,~ \mbfE_{\L}(N(\phi))>1.~\end{eqnarray}
(Note that for a BRW, $\kappa_1(x), \Phi_1(x):=e^{\kappa_1(x)}, \mbfE_{\L}(N(\phi))$ are all constants since there is only one element in the state space of $\L$.) Hence \eqref{Jafassume} is the degenerate version of the Condition 1 and Condition 3 in the present paper.

Furthermore, if the random environment is degenerate, Condition 2 and Condition 4 will hold automatically under Condition 1 and Condition 3. Let us check them one by one. The first equality in \eqref{Jafassume} means the left hand side of \eqref{c2'a} is $0$ thus \eqref{c2'a} holds. We see $\Phi_1(0)\in(0,+\infty)$ since $\Phi_1(0)=\mbfE_{\L}(N(\phi))$ and $\mbfE_{\L}(N(\phi)^{1+\lambda_{14}})<+\infty.$ As $\Phi_1$ is the log-Laplace transform of a measure on $\bfR$ and the condition that $\max\{\Phi_1(0),\Phi_1(\vartheta+\lambda_{13})\}<+\infty,$ we see that $\frac{d^n\Phi(\theta)}{d\theta}$ exists and is finite for any $n\in\bfN$ and $\theta\in(0,\vartheta+\lambda_{13})$, which means that $\eqref{c21a}$ holds and thus \eqref{c2'b} holds.

At last, we show that Condition 4 also holds when the random environment is degenerate. The statement that ``Condition 4 is not true" is equivalent to say ``$\mbfE_{\L}\left(\sum_{i=1}^{N(\phi)}1_{\{\vartheta \zeta_i(\phi)+ \kappa_1(\vartheta)\leq 0\}}\right)=0$". Combining with the first equality in \eqref{Jafassume}, we get $\mbfP_{\L}(\min_{i\leq N(\phi)}\zeta_i(\phi)>-\vartheta^{-1}\kappa_1(\vartheta))=1$ from $\mbfE_{\L}\left(\sum_{i=1}^{N(\phi)}1_{\{\vartheta \zeta_i(\phi)+ \kappa_1(\vartheta)\leq 0\}}\right)=0.$ However, we can see
$\mbfP_{\L}(\min_{i\leq N(\phi)}\zeta_i(\phi)>-\vartheta^{-1}\kappa_1(\vartheta))<1$ since $-\kappa'_1(\vartheta)=\frac{\mbfE_{\L}\left(\sum_{i=1}^{N(\phi)}\zeta_i(\phi)e^{-\vartheta \zeta_i(\phi)}\right)}{\mbfE_{\L}\left(\sum_{i=1}^{N(\phi)}e^{-\vartheta \zeta_i(\phi)}\right)},$ which leads to a contradiction hence we prove the truth of Condition 4.
\end{remark}

\begin{remark}
The impact from the random environment to the survival probability and the extinction rate is reflected by the quantity $\sigma_{A}.$
We should note that all the conclusions in Theorem \ref{1} and Theorem \ref{2} reflect that the event $\{Y_n\geq 1\}$ and $\mathcal{S}$ will happen with smaller probability when $\sigma_{A}$ takes a larger value. As an example, here we compare the critical coefficient $a_c$ in the present paper with the corresponding one (i.e. the critical coefficient of the $\frac{1}{3}$ order of the barrier for BRW, we denote it $a_0$) in Jaffuel \cite{BJ2012}. With the notations in the present paper, the value of $a_0$ is $\frac{3\sqrt[3]{3\pi^2\vartheta^2\kappa''_1(\vartheta)}}{2\vartheta}.$
On the other hand, according to the relationship $\gamma(\beta)\geq \frac{\pi^2(1+\beta^2)}{2}$ which has been shown in \cite{LY201801}, we have $\gamma_{\sigma}\geq \frac{\pi^2(\sigma^2_A+\sigma^2_Q)}{2}$ hence $a_c=\frac{3\sqrt[3]{6\gamma_\sigma}}{2\vartheta}\geq \frac{3\sqrt[3]{3\pi^2\sigma^2_A+3\pi^2\sigma^2_Q}}{2\vartheta}.$ Note that $\sigma^2_Q:=\vartheta^2\mathbf{E}(\kappa''_1(\vartheta)).$ Hence the term containing $\sigma^2_A$ can be seen as an extra increment of the critical coefficient brought by the random environment.

\end{remark}

\begin{remark}
{\rm (1)} The upper bounds in Theorem \ref{2} only need Conditions 1-2. In other words, under Conditions 1-2, it is true that 
$\varlimsup\limits_{n\rightarrow\infty}\frac{\log\mbfP_{\L}(Y_n>0)}{\sqrt[3]{n}}\leq -\vartheta\tilde{q}(0)$ in {\rm(2a)},
$\varlimsup\limits_{n\rightarrow\infty}\frac{\log\mbfP_{\L}(Y_n>0)}{\sqrt[3]{n}}\leq-\sqrt[3]{3\gamma_{\sigma}}$ in {\rm(2b)} and
$\varlimsup\limits_{n\rightarrow\infty}\frac{1}{n^{1/3}}\log P_{\L}\left(n,b\right)\leq -x_b$ in {\rm(2c)}.

{\rm (2)} Set $\varphi_{\L}(i):=-\vartheta^{-1}K_i+\psi(i),$ from the proof of {\rm(2b)} we can see that
~$\eqref{2b}$ still holds if the function $\psi$ satisfies that $\lim\limits_{n\rightarrow +\infty}\frac{\max_{i\leq n}\psi(i)}{n^{1/3}}=0$ and $\psi(i)\geq 0, \forall i\in\bfN$.

{\rm (3)} We remind that the Propositions 3.2-3.6 in Jaffuel \cite{BJ2012} have shown that the integral equation \eqref{2a1} has an unique solution satisfying the boundary condition \eqref{2a2}.
\end{remark}
\begin{remark}\label{5--}
If $N(\phi)$ and $(\zeta_i(\phi), i\leq N(\phi))$ are independent of each other (i.e. the branching law and the displacement law are independent), we do not need the restriction $\mbfE(|\kappa_1(\vartheta+\lambda_4)|^{\lambda_3})+\mbfE(|\kappa_1(\vartheta)|^{\lambda_3})<+\infty$ in Condition 3. The proof of this remark is given in Section 5.
\end{remark}
\begin{remark}\label{5}
We can see that Condition 4 is almost a necessary condition in some extent from the following analysis.
Now we give a new condition which is slightly weaker than Condition 4. The new condition is that there exists a constant $\lambda_5>2, c>0$ such that
\begin{eqnarray}\label{c5-}\mathbf{E}\left(\left|\log \mbfE_{\L}\left(\sum_{i=1}^{N(\phi)}1_{\{\vartheta \zeta_i(\phi)+ \kappa_1(\vartheta)\leq c\}} \right)\right|^{\lambda_5}\right)<+\infty.\end{eqnarray}
(We should note that under the second assumption in Condition 3, Condition 4 implies \eqref{c5-}.) Then we can see that \eqref{c5-}
is a necessary condition for some conclusions in Theorem \ref{1} when the distribution of $\L_1$ is supported in only finite elements (denoted by $\mathfrak{m}_1,\cdots,\mathfrak{m}_k$) in $\Pi.$ 

The explanation is as follows. Note that for the finite environments case, the statement that \eqref{c5-} is not satisfied is equivalent to say that there exists a $j\leq k$ such that $$\mbfE_{\L}\left(\sum_{i=1}^{N(\phi)}1_{\{\vartheta \zeta_i(\phi)+ \kappa_1(\vartheta)\leq c\}}\Big|\L_1=\mathfrak{m}_j\right)=0,$$ which means $\mbfP_{\L}\left(\min_{i\leq N(\phi)}\zeta_i(\phi)>-\vartheta^{-1}K_1+\vartheta^{-1}c|\L_1=\mathfrak{m}_j\right)=1$. Then for the barrier function $\varphi_{\L}(i):=-\vartheta^{-1}K_i+a i^{\alpha}, a>0, \alpha<1$ and an integer $l$ large enough such that $\vartheta^{-1}cl>a l^{\alpha},$ we have \begin{eqnarray}\label{prop4}\mbfP_{\L}(\mathcal{S})=0~~ \text{if}~~ \L_i=\mathfrak{m}_j, ~~\forall i\leq l.\end{eqnarray} But $\mbfP(\{\L_i=\mathfrak{m}_j, \forall i\leq l\})>0$ since we assume that $\L_1$ only takes finite states with positive probability. Hence \eqref{prop4} contradicts Theorem \ref{1} (b) and the $\alpha\in(\frac{1}{3},1)$ part of Theorem \ref{1} (a).
\end{remark}

To prove the Theorem \ref{1} and Theorem \ref{2}, we still need two important lemmas: the many-to-one formula of the time-inhomogeneous bivariate version and the small deviation principle of a random walk with time-inhomogeneous random environment. We will introduce them in Section 3 and Section 4 respectively.

\section{The many-to-one formula of time-inhomogeneous bivariate version}

The many-to-one formula (a kind of measure transformation named from changing all the paths in the random genealogical tree to a random walk) is an essential tool in the study of the branching random walks.  It can be traced down to the early works of Peyri\`{e}re \cite{P1974}~and Kahane and Peyri\`{e}re \cite{KP1976}.
We refer also to \cite{BK2004} for more variations of this result. The version of time-inhomogeneous many-to-one formula has been first introduced in Mallein \cite{M2015a}. On the other hand, for the time-homogeneous case the bivariate version of many-to-one formula can be found in \cite{GHS2011}.  In this paper we need to establish a bivariate version of many-to-one formula in a time-inhomogeneous random environment.
Let $\tau_{n,\L}$ be a random probability measure on $\bfR\times\bfN$ such that for any $x\in\bfR, A\in\bfN,$ we have
\beqlb\label{m1}\tau_{n,\L}\left((-\infty,x]\times[0,A]\right)=\frac{\mbfE_{\L}\left(1_{\{N(u)\leq A\}}\sum^{N(u)}_{i=1}1_{\{\zeta_i(u)\leq x\}}e^{-\vartheta \zeta_i(u)}\right)}{\mbfE_{\L}\left(\sum^{N(u)}_{i=1}e^{-\vartheta \zeta_i(u)}\right)},~~|u|=n-1,\eeqlb
where $\vartheta$ has been introduced in Condition 1.
Hence we can see that the randomness of $\tau_{n,\L}$ comes entirely from $\L_n.$ Moreover, since $N(u)$ only takes values on $\bfN,$ we have
$$\tau_{n,\L}(\bfR\times([0,+\infty)\setminus \bfN))=0, ~~~\rm{\mathbf{P}-a.s.}$$
Under the quenched law $\mbfP_{\L},$ we introduce a series of independent two-dimensional random vectors $\{X_n, \xi_n\}_{n\in\bfN^+}$ whose distributions are $\{\tau_{n,\L}\}_{n\in\bfN^+}.$ Define \beqlb\label{shift0}S_0:=0, ~S_n:=\sum_{i=1}^{n}X_i, ~\forall n\in\bfN^+.\eeqlb Before we give the many-to-one formula to show the relationship between $\{(S_n, \xi_n), n\in\bfN^+\}$ and the BRWre, we still need the shift operator.

 Since the typical realization of $\L$ is a time-inhomogeneous environment, it is necessary to introduce the shift operator~$\mathfrak{T}.$~Define
$$\mathfrak{T}_0\L:=\L,~~\mathfrak{T}\L:=(\L_2,\L_3,\ldots),~~\mathfrak{T}_k:=\mathfrak{T}^{*k},~~\forall k\in\bfN^+,$$
hence~$\mathfrak{T}_k\L=(\L_{k+1},\L_{k+2},\ldots).$ We use~$\mbfP^k_{\L}$ to present the distribution of~$(\mathbf{T},V,\mbfP_{\mathfrak{T}_k\L}).$ Denote $\mbfE^k_{\L}$ the corresponding expectation of~$\mbfP^k_{\L}.$ (Obviously, $\mbfP^0_{\L}$ and $\mbfE^0_{\L}$ are $\mbfP_{\L}$ and $\mbfE_{\L}$ respectively.) Slightly abusing notation we denote
\beqlb\label{shift} S_n:=\sum_{i=1}^{n}X_{k+i}~~{\rm and}~~ \{\xi_{n}\}_{n\in\bfN}:=\{\xi_{k+n}\}_{n\in\bfN}~~{\rm under}~~ \mbfP^k_{\L}.\eeqlb
That is to say, in this scenario and \eqref{m1}, we have ${\rm\mathbf{P}-a.s.},$
\beqlb\label{m2}\mbfP^k_{\L}(X_{1}\leq x, \xi_{1}\leq A)&=&\mbfP_{\L}(X_{k+1}\leq x, \xi_{k+1}\leq A)\nonumber\\&=&\mbfE_{\L}\left(1_{\{N(u)\leq A,|u|=k\}}\sum^{N(u)}_{i=1}1_{\{\zeta_i(u)\leq x\}}e^{-\vartheta\zeta_i(u)-\kappa_{k+1}(\vartheta)}\right)\nonumber
\\&=&\mbfE^{k}_{\L}\left(\1_{\{N(\phi)\leq A\}}\sum^{N(\phi)}_{i=1}\1_{\{\zeta_i(\phi)\leq x\}}e^{-\vartheta \zeta_i(\phi)-\kappa_{k+1}(\vartheta)}\right).\eeqlb
The following formula gives the relationship between $\{(S_n, \xi_n), n\in\bfN^+\}$ and the BRWre.
\begin{lem}\label{mto}
{\bf (Many-to-one)} For any~$n\in\bfN^+,~k\in\bfN,$ a positive sequence~$\{A_i\}_{i\in \bfN^+}$ and a measurable function~$f:\bfR^n\rightarrow [0, +\infty),$ in the
~sense of~~${\rm\mathbf{P}-a.s.},$ we have
\begin{eqnarray}\label{mto0}
&&\mbfE^k_{\L}\left[\sum_{|u|=n}f(V(u_i),1\leq i\leq n)\1_{\{N(u_{i-1})\leq A_i,1\leq i\leq n\}}\right]\nonumber
\\&&~~~~~~~~~~~~~~~~=\mbfE^k_{\L}\left[e^{\vartheta S_n+\sum^{n}_{i=1}\kappa_{k+i}(\vartheta)}f(S_i,1\leq i\leq n)\1_{\{\xi_{i}\leq A_i,1\leq i\leq n\}}\right].\end{eqnarray}
\end{lem}
Borrowing the idea from the proof of \cite[Theorem1.1]{S2015}, we also prove it by induction.

\noindent{\bf Proof of Lemma \ref{mto}}

We prove it by induction on~$n$. For~$n=1$, if $f$ has the form $f(x)=e^{-x}\1_{\{x\leq A\}}$, \eqref{mto0} can be deduced easily by \eqref{m1} and \eqref{m2}. Hence \eqref{mto0} also holds when $f$ is a non-negative measurable function by the standard method.

We assume that \eqref{mto0} holds as~$n=j$.

Now we consider~$n=j+1.$  We should note that the ancestor $\phi$ of the particle system~under $\mbfP^k_{\L}$ can also be viewed as a particle alive at the $k$-th generation under $\mbfP_{\L}.$ For any non-negative measurable function $f$ defined on $\bfR^{j+1},$ we denote $g_x: \bfR^j\rightarrow [0, +\infty)$ by $$g_x(y_1,y_2,\ldots,y_j):=f(x,x+y_1,x+y_2,\ldots,x+y_j).$$ Recalling $X(\phi):=(N(\phi), \zeta_1(\phi),\zeta_2(\phi),\ldots)$, we have
\begin{eqnarray}\label{mto22}
&&\mbfE^k_{\L}\left[\sum_{|u|=j+1}f(V(u_i),1\leq i\leq j+1)\1_{\{N(u_{i-1})\leq A_i,1\leq i\leq j+1\}}\right]\no
\\&=&\mbfE^k_{\L}\left[\mbfE^k_{\L}\left(\sum_{|u|=j+1}f(V(u_i),1\leq i\leq j+1)\1_{\{N(u_{i-1})\leq A_i,1\leq i\leq j+1\}}\Bigg|X(\phi)\right)\right]\no
\\&=&\mbfE^k_{\L}\left[\1_{\{N(\phi)\leq A_1\}}\sum_{|v|=1}\mbfE^{^{k+1}}_{\L}\left(\sum_{|u|=j}g_{_{V(v)}}(V(u_i)-V(v),1\leq i\leq j)\1_{\{N(u_{i-1})\leq A_{i+1},1\leq i\leq j\}}\right)\right]\no 
\\&=&\mbfE^k_{\L}\left[\1_{\{N(\phi)\leq A_1\}}\sum_{|v|=1}\mbfE^{^{k+1}}_{\L}\left(e^{\vartheta S_j+\sum\limits^{j}_{i=1}\kappa_{_{k+1+i}}(\vartheta)}g_{_{V(v)}}(S_i,1\leq i\leq j)\1_{\{\xi_i\leq A_{i+1},1\leq i\leq j\}}\right)\right]\no
\\&=&\mbfE^k_{\L}\left[\1_{\{\xi_1\leq A_1\}}e^{\vartheta S_{1}+\kappa_{_{k+1}}(\vartheta)}\mbfE^{^{k+1}}_{\L}\left(e^{\vartheta S_j+\sum\limits^{j+1}_{i=2}\kappa_{_{k+i}}(\vartheta)}g_{_{S_{1|k}}}(S_i,1\leq i\leq j)\1_{\{\xi_i\leq A_{i+1},1\leq i\leq j\}}\right)\right]\no
\\&=&\mbfE^k_{\L}\Bigg[\1_{\{\xi_1\leq A_1\}}e^{\vartheta S_{1}+\kappa_{_{k+1}}(\vartheta)}\no
\\&~~&~~~~~~~~~~~~~\times\mbfE^{k}_{\L}\Bigg(e^{\vartheta (S_{j+1}-S_1)+\sum\limits^{j+1}_{i=2}\kappa_{_{k+i}}}f(S_i,1\leq i\leq j+1)\1_{\{\xi_i\leq A_{i},2\leq i\leq j+1\}}\Bigg|S_1,\xi_1\Bigg)\Bigg]\no
\\&=&\mbfE^k_{\L}\left[e^{\vartheta S_{1}+\kappa_{k+1}(\vartheta)}e^{\vartheta (S_{j+1}-S_1)+\sum^{j}_{i=1}\kappa_{k+1+i}}f(S_i,1\leq i\leq j+1)\1_{\{\xi_i\leq A_{i},1\leq i\leq j+1\}}\right].
\end{eqnarray}
It should be noted that the notation $S_{1|k}$ in $\mbfE^{k+1}_{\L}(\cdot)$ on the fifth line of \eqref{mto22} is to emphasize that this is the $S_1$ under $\mbfE^{k}_{\L}.$
 For the above equalities, the first one and the last one are due to the smoothness of the conditional expectation. The second one and the penultimate one are obtained by the Markov property of BRWre and $\{S_i\}_{i\in\bfN}$ respectively. Let the $n$ in Lemma \ref{mto} be $j$ (resp. $n$ be $1$), we can get the third one (resp. the fourth one) by using \eqref{mto0}.\qed
\section{The small deviation principle for the associated random walk with random environment in time (RWre)}

Lv and Hong \cite{Lv201802} considered the small deviation principle for RWre. In this section, we give some corollaries of the main results in Lv and Hong \cite{Lv201802}, which will be the essential tools for the barrier problem of BRWre. For the case of time-homogeneous, it is a standard and effective way to solve the barrier problem of BRW by applying some estimates on the so-called associated random walk. Let us first give a short review on it.

The small deviation principle for the random walk $\{V_n\}_{n\in\bfN}$ with i.i.d. random increments has been given in Mogul'ski\u{\i} \cite{Mog1974}. When $V_{i+1}-V_{i}$ has expectation $0$ and variance $\sigma^2<+\infty$, Mogul'ski\u{\i} showed that for two continuous functions $g, h$ such that $g(0)<0<h(0), g(s)<h(s), \forall s\in[0,1],~\alpha\in(0,\frac{1}{2}), x\in\bfR,$
\beqlb\label{1.1Mo}\lim\limits_{n\rightarrow +\infty}\frac{\log\mbfP(\forall_{i\leq n}V_i\in[g(i/n)n^{\alpha},h(i/n)n^{\alpha}]|V_0=x)}{n^{1-2\alpha}}=\frac{-\pi^{2}\sigma^2}{2}\int_0^1\frac{1}{(h(s)-g(s))^2}ds.\eeqlb
(The small deviation principle focuses on the probability that a stochastic process has fluctuations below its natural scale. Therefore, we call \eqref{1.1Mo} a small deviation principle according to the setting $\alpha\in(0,\frac{1}{2})$ and the central limit theorem.)

 By the time-homogeneous many-to-one formula (see \cite[Theorem 1.1]{S2015}), we can construct a bijection between a BRW and a random walk with i.i.d. increments (which is usually called the associated random walk).
Based on this relationship, Jaffuel \cite{BJ2012} and A\"{\i}d\'{e}kon, Jaffuel \cite{AJ2011} have studied the barrier problems of BRW by applying \eqref{1.1Mo}. 

In the present paper we consider the barrier problem of BRWre. We use the many-to-one formula given in Lemma \ref{mto} to construct a corresponding associated random walk for BRWre. In the next section, we will show that the corresponding associated random walk for BRWre is just the RWre studied in Lv and Hong \cite{Lv201802}.

\subsection{The associated RWre and its small deviation principle}

Let us give the definition of RWre. We denote $\mu:=\{\mu_n\}_{n\in\bfN^+}$ an i.i.d. sequence with values in the space of probability measures on $\bfR.$ Conditioned to a realization 
of $\mu,$ we sample $\{\mathcal{V}_n\}_{n\in\bfN^+}$ a sequence of independent random variables such that for every $n\in\bfN^+,$ the law of $\mathcal{V}_n$ is the realization of $\mu_n.$
Set \beqlb\label{sec-sn} \mathcal{V}_0=x\in\bfR,~~\tilde{V}_n:=\mathcal{V}_0+\sum_{i=1}^{n}\mathcal{V}_i.\eeqlb We call $\{\tilde{V}_n\}_{n\in\bfN}$ the {\it random walk with time-inhomogeneous random environment $\mu$}, which is often abbreviated as RWre in the rest of this paper\footnote{Note that the process we consider in the present paper is not the classical
 random walk in random environment which has been well-studied in Zeitouni \cite{Z2004} and many other papers. 
 For the classical random walks in random environment, the random environment is either purely spatial or space-time. However, in our model, the random environment is in time.}. 

 Note that the $\{S_n\}_{n\in\bfN}$ defined in \eqref{mto0} is a RWre with
time-inhomogeneous random environment $\L$ (More precisely, the random environment is the $\{\tau_{n,\L}\}_{n\in\bfN^+}$ with $A=+\infty$ which is defined in \eqref{m1}. Note that $\{\tau_{n,\L}\}_{n\in\bfN^+}$ is totally determined by $\L$ hence we say ``with random environment $\L$").
That is to say, Lemma \ref{mto} constructs the relationship between BRWre and RWre. The studying of small deviation principle for RWre is an important step to solve the barrier problem of BRWre.

Now we introduce the associated RWre. Recall \eqref{shift0} and \eqref{shift} and define \begin{eqnarray}\label{mtoT}
T_n:=\vartheta S_n+K_{n+k}-K_{k} ~~{\rm under}~~ \mbfP^k_{\L},
\end{eqnarray} where $K_n,$ $\vartheta$ and $S_n$ have been defined in \eqref{sect-2}, \eqref{Kn} and Lemma \ref{mto}. Obviously, it is a RWre (with random environment $\mathfrak{T}_k\L$). Next we will show under Conditions 1-2 in the present paper, the $\{T_n\}_{n\in\bfN}$ defined in \eqref{mtoT} satisfies all the basic assumptions in \cite{Lv201802} (thus we can apply the main result in \cite{Lv201802} to $\{T_n\}_{n\in\bfN}$).
 \begin{prop}\label{mt1}~
{\rm (i)}~If~the BRWre satisfies Condition 1, then the associated RWre $T$ satisfies $\mathbf{E}[(T_1-\mbfE_{\L}T_1)^2]>0$ and $\mbfE(T_1)=0$.

{\rm (ii)}~If~the BRWre satisfies Condition 2, then the associated RWre $T$ satisfies $\mbfE([\mbfE_{\L}(T_1)]^{\lambda_1})<+\infty,$ and $\mathbf{E}\left\{[(T_1-\mbfE_{\L}T_1)^{\lambda_2}]^{\lambda_1}\right\}<+\infty$.

In conclusion, if~the BRWre satisfies Conditions 1-2, then the associated RWre $T$ satisfies \cite[Assumptions (H1)-(H3)]{Lv201802}.
\end{prop}

\noindent{\bf Proof of Proposition \ref{mt1}}

{\rm (i)} Recall that $T_n:=\vartheta S_n+K_n, X_n:=S_n-S_{n-1}, K_n:=\sum^n_{i=1}\kappa_i(\vartheta), \kappa(\vartheta):=\mbfE(\kappa_1(\vartheta))$ and $\{\kappa_n(\vartheta)\}_{n\in\bfN}$ is a sequence of i.i.d. random variables. By the many-to-one formula (Lemma \ref{mto}), we see
\begin{eqnarray}\label{mto1}\mbfE_{\L}(X_n)=\frac{\mbfE_{\L}\left(\sum_{i=1}^{N(u)}\zeta_i(u)e^{-\vartheta \zeta_i(u)}\right)}{\mbfE_{\L}\left(\sum_{i=1}^{N(u)}e^{-\vartheta \zeta_i(u)}\right)}=-\kappa'_n(\vartheta),~~|u|=n-1.\end{eqnarray}
Then it is not hard to see $\mathbf{E}(T_1)=0$ since
\begin{eqnarray}\label{mto2}
\mathbf{E}(T_n-T_{n-1})&=&\vartheta\mathbf{E}(X_n)+\mathbf{E}(K_n-K_{n-1})=\vartheta\mathbf{E}(\mbfE_{\L}(X_n))+\mathbf{E}(\kappa_n(\vartheta))\no%
\\&=&-\vartheta\mathbf{E}(\kappa'_n(\vartheta))+\mathbf{E}(\kappa_n(\vartheta))
=-\vartheta\kappa'(\vartheta)+\kappa(\vartheta)=0.
\end{eqnarray}
Lemma \ref{mto} also tells that
\begin{eqnarray}\label{erjiedao}
\kappa''_1(\vartheta)
&=&\frac{\mbfE_{\L}\left(\sum_{i=1}^{N(\phi)}\zeta^2_i(\phi)e^{-\vartheta \zeta_i(\phi)}\right)\mbfE_{\L}\left(\sum_{i=1}^{N(\phi)}e^{-\vartheta \zeta_i(\phi)}\right)-\left[\mbfE_{\L}\left(\sum_{i=1}^{N(\phi)}\zeta_i(\phi)e^{-\vartheta \zeta_i(\phi)}\right)\right]^2}{\left[\mbfE_{\L}\left(\sum_{i=1}^{N(\phi)}e^{-\vartheta \zeta_i(\phi)}\right)\right]^2}\no
\\&=&\mbfE_{\L}(S^2_1)-(\mbfE_{\L}S_1)^2.\end{eqnarray}
Then it is true that \begin{eqnarray}\label{qq}
\mbfE_{\L}((T_1-\mbfE_{\L}T_1)^2)=\mbfE_{\L}((\vartheta S_1-\vartheta\mbfE_{\L}S_1)^2)=\vartheta^2\kappa''_1(\vartheta).\end{eqnarray}
Hence we get $\mbfE((T_1-\mbfE_{\L}T_1)^2)=\vartheta^2\kappa''(\vartheta)>0$ by the Proposition \ref{propcon1}.


{\rm (ii)} By Lemma \ref{mto} we also have \begin{eqnarray}\label{et}
\mbfE_{\L}T_1=\kappa_1(\vartheta)-\vartheta\kappa'_1(\vartheta).
\end{eqnarray} Hence $\mbfE([\mbfE_{\L}T_1]^{\lambda_1})<+\infty$ is equivalent to \eqref{c2'a}. Moreover, according to the many-to-one formula we can see directly that $\mathbf{E}\left\{[(T_1-\mbfE_{\L}T_1)^{\lambda_2}]^{\lambda_1}\right\}<+\infty$  is equivalent to \eqref{c2'b}.\qed

We stress that by the definitions in \eqref{AQ} and the above proof, one can see that
$$\sigma^2_A=\mbfE[(\mbfE_{\L}T_1)^{2}]<+\infty,~~ \sigma^2_Q=\mathbf{E}\left[(T_1-\mbfE_{\L}T_1)^2\right]\in(0,+\infty).$$  Recalling the notation  $\gamma_{\sigma}:=\sigma^2_Q\gamma\left(\frac{\sigma_A}{\sigma_Q}\right)$ and applying the main result in \cite{Lv201802} to $\{T_n\}_{n\in\bfN}$ we get the following result.



\begin{thm}\label{mogc1}
(Lv and Hong \cite[Corollary 2]{Lv201802}) 
~Let $g(s), h(s)$ be two continuous functions on $[0,1]$ and $g(s)<h(s)$ for any $s\in [0,1].$ We set $g(0)< a_0\leq b_0 <h(0),~ g(1)\leq a'<b'\leq h(1)$
~and~ $C_{g,h}:=\int_{0}^{1}\frac{1}{[h(s)-g(s)]^2}ds.$ Let $\{t_n\}_{n\in\bfN}$ be a sequence of non-negative integers and $\bar{t}_n:=t_n+n$. 
For any $\alpha\in (\frac{1}{\lambda_1},\frac{1}{2}),$ 
(for the present paper we only need the case $\alpha=\frac{1}{3}$ hence we require $\lambda_1>3$ in Condition 2,) under Conditions 1-2 we have
\beqlb\label{mogc1u}\lim\limits_{n\rightarrow +\infty}\frac{\sup\limits_{x\in\bfR}\log\mbfP_\L
\Big(\forall_{0\leq i\leq n}T_{t_n+i}\in \big[g\big(\frac{i}{n}\big)n^\alpha,h\big(\frac{i}{n}\big)n^\alpha\big]\Big|T_{t_n}=x\Big)}{n^{1-2\alpha}}= -C_{g,h}\gamma_{\sigma},~{\rm \mathbf{P}-a.s.,}\eeqlb
\beqlb\label{mogc1l}&&\lim\limits_{n\rightarrow +\infty}\frac{\inf\limits_{x\in[a_0 n^{\alpha}, b_0 n^{\alpha}]}\log \mbfP_\L
\left(\forall_{0\leq i\leq n} T_{t_n+i}\in \left[g\left(\frac{i}{n}\right)n^\alpha,h\left(\frac{i}{n}\right)n^\alpha\right], T_{\bar{t}_n}\in [a'n^\alpha,b'n^\alpha] \Bigg|T_{t_n}=x\right)}{n^{1-2\alpha}}\no
\\&&~~~~~~~~~~~~~~~~~~~~~~~~~~~~~~~~~~~~~~~~~~= -C_{g,h}\gamma_{\sigma},~~~~{\rm \mathbf{P}-a.s.},\eeqlb
where $\{T_n\}_{n\in\bfN}$ is the one in \eqref{mtoT}.
\end{thm}

Therefore, if the random environment $\L$ is degenerate (thus $\sigma_A=0$), then we can see that Theorem \ref{mogc1} is consistent with the Mogul'ski\u{\i} estimate \eqref{1.1Mo} since $\gamma(0)=\frac{\pi^2}{2}$ (see \eqref{bmt}).

However, Theorem \ref{mogc1} still can not be applied directly to prove the main results in the present paper. Hence we need the following three useful corollaries of Theorem \ref{mogc1}. 
\subsection{Some Corollaries of Theorem \ref{mogc1}}

In the forthcoming three corollaries, we will see why we need Condition 3 and Condition 4. In short, Condition 3 allows us to add some extra events in the lower bound \eqref{mogc1l} and Condition 4 ensures that \eqref{mogc1l} holds even though $b_0=h(0).$

From now on, we always set $\alpha=\frac{1}{3}$ since in the present paper we only need the case $\alpha=\frac{1}{3}.$ But we point out that Corollaries \ref{mogc01}-\ref{mogc2} (resp. Corollary \ref{mog0}) will also hold for $\alpha\in(\frac{1}{\lambda_1},\frac{1}{2})$ (resp. $\alpha\in(\frac{1}{\lambda_1},\frac{1}{3}]$), where $\lambda_1$ is the one defined in Condition 2. 

\begin{cor}\label{mogc01}
$\{\xi_n\}_{n\in\bfN}$ has been introduced in Section 3. 
The setting of $g,~h,~ a_0,b_0, a',b', t_n, \bar{t}_n$ and~$C_{g,h}$ are the same as what we introduce in Theorem \ref{mogc1}. Let $v\in(\frac{2}{\lambda_3},+\infty),$ where $\lambda_3$ is the one introduced in Condition 3.
Then under Conditions 1-3 we have
\beqlb\label{mogc1l+}&&\lim\limits_{n\rightarrow +\infty}\frac{\inf\limits_{x\in[a_0 n^{\alpha}, b_0 n^{\alpha}]}\log \mbfP_\L
\left( \substack{\forall_{0\leq i\leq n}~~\frac{T_{t_n+i}}{n^\alpha}\in \left[g\left(\frac{i}{n}\right),h\left(\frac{i}{n}\right)\right],\\ \frac{T_{\bar{t}_n}}{n^\alpha}\in [a',b'],~~~ \xi_{t_n+i}\leq \exp\{n^{v}\}~~} \Bigg|T_{t_n}=x\right)}{n^{1-2\alpha}}\no
\\&&~~~~~~~~~~~~~~~~~~~~~~~~~~~~~~~~~~~~~~~~~~= -C_{g,h}\gamma_{\sigma},~~~~{\rm \mathbf{P}-a.s.}\eeqlb
\end{cor}

\begin{cor}\label{mogc2}
Let $l,m,N\in\bfN,$ $0\leq l< m\leq N, v\in(\frac{2}{\lambda_3}, +\infty).$
Let $g(s), h(s)$ be two continuous functions on $[0,1]$ such that $g(s)<h(s), \forall s\in [0,1],~ g(l/N)< a_0\leq b_0 <h(l/N), ~g(m/N)\leq a'<b'\leq h(m/N).$ 
For $0\leq z_1<z_2\leq 1$,~define $C^{z_1,z_2}_{~g,~h}:=\int_{z_1}^{z_2}\frac{1}{[h(s)-g(s)]^2}ds.$ 
Under Conditions 1-3, we have
$$\lim\limits_{k\rightarrow +\infty}\frac{\sup\limits_{x\in\bfR}\log\mbfP_\L
\left(\forall_{lk\leq i\leq mk}\frac{T_{i}}{(Nk)^\alpha}\in [g(\frac{i}{Nk}),h(\frac{i}{Nk})]\big|T_{lk}=x\right)}{(Nk)^{1-2\alpha}}= -C^{\frac{l}{N},\frac{m}{N}}_{~g,~h}\gamma_{\sigma},~~~{\rm \mathbf{P}-a.s.}$$
\beqnn&&\lim\limits_{k\rightarrow +\infty}\frac{\inf\limits_{x\in[a_0 (Nk)^{\alpha}, b_0 (Nk)^{\alpha}]}\log \mbfP_\L
\left(\substack{\forall_{lk\leq i\leq mk}~\frac{T_{i}}{(Nk)^\alpha}\in [g(\frac{i}{Nk}),h(\frac{i}{Nk})],\\~~\frac{T_{mk}}{(Nk)^\alpha}\in [a',b'],~~~\xi_i\leq \exp\{(Nk)^{v}\}}\Bigg|T_{lk}=x\right)}{(Nk)^{1-2\alpha}}
= -C^{\frac{l}{N},\frac{m}{N}}_{~g,~h}\gamma_{\sigma},~{\rm \mathbf{P}-a.s.}
\eeqnn
\end{cor}

We should note that the above two corollaries both need the condition $g(0)< a_0\leq b_0 <h(0).$ In the next corollary we consider the case $h(0)=b_0.$  The following lemma is a necessary preparation for the next corollary.
\begin{lemma}\label{lem4}
Let $\{X_n\}_{n\in\bfN}$ be a sequence of i.i.d. random variables and $\{t_n\}_{n\in\bfN}$ a sequence of non-negative integers. Assume that there exists an $\epsilon>0$ such that $\bfE(|X_1|^{2+\epsilon})<+\infty,$ then \beqlb\label{l4l}\lim\limits_{n\rightarrow+\infty}n^{-1}\sum_{i=t_n+1}^{t_n+n}X_i=\bfE(X_1),~{\rm a.s.}\eeqlb
\end{lemma}
\begin{cor}\label{mog0}
The setting of $g,~h,~ a',b', t_n$ and $C_{g,h}$ are the same as what we introduce in Theorem 4.1.
Suppose that $h(s)\geq h(0), \forall s\in(0,1].$
 Under the Conditions 1-4, we have
\begin{eqnarray}\label{mog.1}\lim \limits_{n\rightarrow +\infty}\frac{\log \mbfP_\L
\left(\substack{\forall_{0\leq i\leq n}~\frac{T_{t_n+i}}{n^\alpha}\in \left[g\left(\frac{i}{n}\right),h\left(\frac{i}{n}\right)\right],~\\ \frac{T_{t_n+n}}{n^\alpha}\in [a',b'],~~~\xi_{t_n+i}\leq \exp\{n^{v}\}}\Bigg|T_{t_n}=h(0)n^{\alpha}\right)}{n^{1-2\alpha}}
= -C_{g,h}\gamma_{\sigma},~~{\rm \mathbf{P}-a.s.}\end{eqnarray}
Especially,
if $t_n\equiv k, k\in\bfN$  
then \eqref{mog.1} still holds under Conditions 1-4 even though the $\lambda_5$ described in Condition 4 only satisfies $\lambda_5\geq 1$.
\end{cor}

\subsection{The proofs of Corollaries \ref{mogc01}-\ref{mog0} and Lemma \ref{lem4}}

\noindent{\bf Proof of Corollary \ref{mogc01}}

The only gap between \eqref{mogc1l} and \eqref{mogc1l+} is the term $``\xi_{t_n+i}\leq e^{n^{v}}".$ We will cross the gap by using Condition 3.
With the help of \eqref{mogc1u}, to prove \eqref{mogc1l+} we only need to show
\begin{eqnarray}\label{mog.1>}\varliminf \limits_{n\rightarrow +\infty}\frac{\inf\limits_{x\in[a_0 n^{\alpha}, b_0 n^{\alpha}]}\log \mbfP_\L
\left(\substack{\forall_{0\leq i\leq n}~\frac{T_{t_n+i}}{n^\alpha}\in \left[g\left(\frac{i}{n}\right),h\left(\frac{i}{n}\right)\right],~\\ \frac{T_{t_n+n}}{n^\alpha}\in [a',b'],~~~\xi_{t_n+i}\leq e^{n^{v}}}\Bigg|T_{t_n}=x\right)}{n^{1-2\alpha}}
\geq -C_{g,h}\gamma_{\sigma},~~{\rm \mathbf{P}-a.s.}\end{eqnarray}
First, we consider the case that $g(x)=a, h(x)=b, \forall x\in[0,1],$ where $a<a_0\leq b_0<b,~ a\leq a'<b'\leq b.$ The proof have the same spirit as the proof of the lower bound in \cite[Theorem 2]{Lv201802}. Now we show the adjustments required in the proof of \cite[Theorem 2]{Lv201802}. 
Choose $a'',b''$ satisfying $a'<a''<b''<b'.$
Let $D\in\bfN^+, ~J:=\lfloor Dn^{2\alpha}\rfloor,~ K:=\left\lfloor \frac{n}{J} \right\rfloor,$ $t_{n,k}:=t_n+kJ.$
By Markov property we have 
\begin{eqnarray}
&&\inf\limits_{x\in[a_0 n^{\alpha}, b_0 n^{\alpha}]}\mbfP_\L
\Big(\forall_{t_n\leq i\leq \bar{t}_n}~ T_{i}\in[an^{\alpha},bn^{\alpha}],~T_{\bar{t}_n}\in[a'n^{\alpha}, b'n^{\alpha}],~\xi_i\leq e^{n^{v}} \Big|T_{t_n}=x\Big)\no
\\&\geq&\prod\limits_{k=0}\limits^{K-1}\inf\limits_{x\in [a_0 n^\alpha,b_0 n^\alpha]}\left[\mbfP_\L\left(\substack{\forall_{i\leq J}~T_{t_{n,k}+i}\in[a n^\alpha, bn^\alpha],\\~T_{t_{n,k+1}}\in[a''n^\alpha,b''n^\alpha]}\Bigg|T_{t_{n,k}}=x\right)-\sum_{i=t_{n,k}+1}^{t_{n,k+1}}\mbfP_{\L}(\xi_i>e^{n^v})
\right]\no
\\&\times&\inf\limits_{x\in [a'' n^\alpha,b'' n^\alpha]} \left[\mbfP_\L\left(\substack{\forall_{i\leq \bar{t}_n-t_{n,K}}~ T_{t_{n,K}+i}\in[an^\alpha, b n^\alpha],\\ T_{\bar{t}_n}\in[a'n^\alpha,b'n^\alpha]} \Bigg|T_{t_{n,K}}=x\right)-\sum_{i=t_{n,K}+1}^{\bar{t}_n}\mbfP_{\L}(\xi_i>e^{n^v})\right].\no
\end{eqnarray}
By the many-to-one formula, we have
$$\mbfP_{\L}(\xi_i>e^{n^v})=\frac{\mbfE_{\L}\left(\1_{\{N(u)>e^{n^v}\}}\sum\limits_{j=1}^{N(u)}e^{-\vartheta\zeta_{j}(u)}\right)}{\mbfE_{\L}\left(\sum\limits_{j=1}^{N(u)}e^{-\vartheta\zeta_{j}(u)}\right)} =\frac{\mbfE_{\L}\left(\1_{\{N(u)>e^{n^v}\}}\sum\limits_{j=1}^{N(u)}e^{-\vartheta\zeta_{j}(u)}\right)}{e^{\kappa_i(\vartheta)}},~~|u|=i-1.$$
Let $v_1:=\frac{\lambda_4}{\vartheta+\lambda_4},$ where the $\lambda_4$ has been introduced in Condition 3. By H\"{o}lder's inequality we get
\beqlb\label{Holder}&&\mbfE_{\L}\left(\1_{\{N(u)>e^{n^v}\}}\sum_{j=1}^{N(u)}e^{-\vartheta\zeta_{j}(u)}\right)\no
\\&=& \mbfE_{\L}\left(\1_{\{N(u)>e^{n^v}\}}N(u)^{v_1}N(u)^{-v_1}\sum_{j=1}^{N(u)}e^{-\vartheta\zeta_{j}(u)}\right)\no
\\&\leq&\left[\mbfE_{\L}\left(\1_{\{N(u)>e^{n^v}\}}N(u)\right)\right]^{v_1}\left[\mbfE_{\L}\left(\left(\sum_{j=1}^{N(u)}N(u)^{-v_1}e^{-\vartheta\zeta_{j}(u)}\right)^{\frac{1}{1-v_1}}\right)\right]^{1-v_1}.
\eeqlb
The fact $\frac{1}{1-v_1}>1$ implies that $$\left(\sum_{j=1}^{N(u)}N(u)^{-v_1}e^{-\vartheta\zeta_{j}(u)}\right)^{\frac{1}{1-v_1}}\leq N(u)^{\frac{1}{1-v_1}-1}\sum_{j=1}^{N(u)}\left(N(u)^{-\frac{v_1}{1-v_1}}e^{-\frac{\vartheta\zeta_{j}(u)}{1-v_1}}\right)=\sum_{j=1}^{N(u)}e^{-(\vartheta+\lambda_4)\zeta_{j}(u)}.$$
Hence by Markov property and the above inequality we get
\beqlb&&\mbfE_{\L}\left(\1_{\{N(u)>e^{n^v}\}}\sum_{i=1}^{N(u)}e^{-\vartheta\zeta_{i}(u)}\right)\no
\\&\leq&\left[\mbfE_{\L}\left(\1_{\{N(u)>e^{n^v}\}}\frac{N(u)^{1+\lambda_3}}{e^{\lambda_3n^v}}\right)\right]^{v_1}(e^{\kappa_i(\vartheta+\lambda_4)})^{1-v_1}\no
\\&\leq&e^{-\lambda_3v_1n^v}\mbfE_{\L}(N(u)^{1+\lambda_3})^{v_1}(e^{\kappa_i(\vartheta+\lambda_4)})^{1-v_1}.\eeqlb
Choose $v_2\in(\frac{2}{\lambda_3},v)$ and define \beqlb\label{tilIn}\hat{I}_n:=\left\{\max_{|u|\leq n}\mbfE_{\L}(N(u)^{1+\lambda_3})\leq e^{n^{v_2}}, \max_{i\leq n}[(1-v_1)\kappa_i(\vartheta+\lambda_4)-\kappa_i(\vartheta)]\leq n^{v_2}\right\}.\eeqlb
Hence on the event $\hat{I}_n,$ for $n$ large enough, we have
\begin{eqnarray}\label{PROP3}
&&\inf\limits_{x\in[a_0 n^{\alpha}, b_0 n^{\alpha}]}\mbfP_\L
\Big(\forall_{t_n\leq i\leq \bar{t}_n}~ T_{i}\in[an^{\alpha},bn^{\alpha}],~T_{\bar{t}_n}\in[a'n^{\alpha}, b'n^{\alpha}],~\xi_i\leq e^{n^{v}} \Big|T_{t_n}=x\Big)\no
\\&\geq&\prod\limits_{k=0}\limits^{K-1}\left[\inf\limits_{x\in [a_0 n^\alpha,b_0 n^\alpha]} \mbfP_\L\left(\substack{\forall_{i\leq J}~T_{t_{n,k}+i}\in[a n^\alpha, bn^\alpha],\\~T_{t_{n,k+1}}\in[a''n^\alpha,b''n^\alpha]}\Bigg|T_{t_{n,k}}=x\right)-e^{-\frac{\lambda_3v_1n^v}{2}}
\right]\no
\\&\times&\inf\limits_{x\in [a'' n^\alpha,b'' n^\alpha]} \left[\mbfP_\L\left(\substack{\forall_{i\leq \bar{t}_n-t_{n,K}}~ T_{t_{n,K}+i}\in[an^\alpha, b n^\alpha],\\ T_{\bar{t}_n}\in[a'n^\alpha,b'n^\alpha]} \Bigg|T_{t_{n,K}}=x\right)-e^{-\frac{\lambda_3v_1n^v}{2}}\right].
\end{eqnarray}
That is to say, \cite[(3.13)]{Lv201802} still holds even though we change the left hand side of \cite[(3.13)]{Lv201802} from
$\mbfP_\mu\left(\substack{\forall_{i\leq \lfloor Dn^{2\alpha}\rfloor}~ S_{t_{n,k}+i}\in[a n^\alpha, bn^\alpha],\\ ~S_{t_{n,k+1}}\in[a''n^\alpha,b''n^\alpha] }\Big|S_{t_{n,k}}=x\right)$
to $\mbfP_\L\left(\substack{\forall_{i\leq J}~ T_{t_{n,k}+i}\in[a n^\alpha, bn^\alpha],\\ ~T_{t_{n,k+1}}\in[a''n^\alpha,b''n^\alpha] ,~\xi_{t_{n,k}+i}\leq e^{n^v}}\Big|T_{t_{n,k}}=x\right)$ (the $\{S_n\}_{n\in\bfN}$ with random environment $\mu$ in \cite{Lv201802} and the $\{T_n\}_{n\in\bfN}$ with random environment $\L$ in the present paper satisfy the same assumptions). Then according to the method of the proof of \cite[Theorem 2]{Lv201802}, we only need to show $\lim\limits_{n\rightarrow+\infty}\1_{\hat{I}_n}=1, ~{\rm \mathbf{P}-a.s.}$

Note that $\{\mbfE_{\L}(N(u)^{1+\lambda_4})\leq e^{n^{v_2}}\}=\{\max\{\mbfE_{\L}(N(u)^{1+\lambda_4}),1\}\leq e^{n^{v_2}}\}$, hence
$$\mathbf{P}(\hat{I}^c_n)\leq n\cdot n^{-\lambda_3v_2}\left[\mbfE([\log^+\mbfE_{\L}(N(u)^{1+\lambda_4})]^{\lambda_3})+\mbfE(|\kappa_1(\vartheta)|^{\lambda_3})+\mbfE(|\kappa_1(\vartheta+\lambda_4)|^{\lambda_3})\right].$$
Note that $\lambda_3v_2>2$, then by Borel-Cantelli 0-1 law we get $\lim\limits_{n\rightarrow+\infty}\1_{\hat{I}_n}=1, ~{\rm \mathbf{P}-a.s.}$  Hence we have $\lim\limits_{n\rightarrow+\infty}\1_{\hat{I}_n\cap H_n}=1, ~{\rm \mathbf{P}-a.s.},$ where the event $H_n$ is defined in \cite[(3.19)]{Lv201802} and we have shown that $\lim\limits_{n\rightarrow+\infty}\1_{H_n}=1, ~{\rm \mathbf{P}-a.s.}$ in \cite{Lv201802}. According to the facts that \cite[(3.13)]{Lv201802} still holds in our context and $\lim\limits_{n\rightarrow+\infty}\1_{\hat{I}_n\cap H_n}=1, ~{\rm \mathbf{P}-a.s.}$,  we can get \eqref{mogc1l+} in the case $g(x)=a, h(x)=b, \forall x\in[0,1]$ by using the method in \cite[Theorem 2]{Lv201802}. Furthermore, by the method used in \cite[Corollary 2]{Lv201802}, we can let $g,h$ be any continuous functions on $[0,1]$ as long as $g(0)< a_0\leq b_0 <h(0),~ g(1)\leq a'<b'\leq h(1).$ \qed


\noindent{\bf Proof of Corollary \ref{mogc2}}

Define two continuous functions on~$[0,1]$ $$~\tilde{g}(x):=\left(\frac{N}{m-l}\right)^\alpha g\left(\left(x+\frac{l}{m-l}\right)\frac{m-l}{N}\right),
~\tilde{h}(x):=\left(\frac{N}{m-l}\right)^\alpha h\left(\left(x+\frac{l}{m-l}\right)\frac{m-l}{N}\right).$$ Then the event
$$\left\{\forall_{lk\leq i\leq mk}~\frac{T_{i}}{(Nk)^\alpha}\in \left[g\left(\frac{i}{Nk}\right),~h\left(\frac{i}{Nk}\right)\right]\right\}$$
can be rewritten as
\begin{eqnarray}\label{3.4.6}\left\{\forall_{i\leq mk-lk}~\frac{T_{i+lk}}{(mk-lk)^\alpha}\in\left[\tilde{g}\left(\frac{i}{mk-lk}\right),\tilde{h}\left(\frac{i}{mk-lk}\right)\right]\right\}.\no\end{eqnarray}
 Then we replace the time length $n$ and the starting time $t_n$ in \eqref{mogc1u} and \eqref{mogc1l+} by $n:=mk-lk$ and $t_n:=\frac{ln}{m-l}=lk$ respectively. 
 Then we can deduce that the limits in \eqref{mogc1u} and~\eqref{mogc1l+} are both 
$C^{0,1}_{\tilde{g},\tilde{h}}\gamma_\sigma\left(\frac{m-l}{N}\right)^{1-2\alpha}$. 
At last, by some standard calculations we get $\left(\frac{m-l}{N}\right)^{1-2\alpha}C^{0,1}_{\tilde{g},\tilde{h}}=C^{\frac{l}{N},\frac{m}{N}}_{~g,~h},$ which completes the proof.\qed 

\noindent{\bf Proof of Lemma \ref{lem4}}

We know \eqref{l4l} can be obtained directly by Borel-Cantelli 0-1 law if $\bfE(X^{4}_1)<+\infty.$ But using the strong approximation method we only need the assumption $\bfE(X^{2+\epsilon}_1)<+\infty.$ Define $\bar{X}_i:=X_i-\bfE(X_1).$ By \cite[Theorem 1]{Sak2006} we can construct a standard Brownian motion $W$ such that
$$\forall x>0, ~n\in\bfN^{+},~\bfP\left(\max_{k\leq n}\left|\sum_{i=1}^{k}\bar{X}_i-W_{k\sigma^2_{X}}\right|\geq C^*(2+\epsilon)x\right)\leq \frac{2n\bfE(\bar{X}_i^{2+\epsilon})}{x^{2+\epsilon}},$$
where $\sigma^2_{X}$ is the variance of $X_1$ and $C^*$ is a positive absolute constant. Moreover, by the Cs\"{o}rg\H{o} and R\'{e}v\'{e}sz's estimation \cite[Lemma 1]{CR1979} we can find two constants $c_s,c'_s>0$ such that
$\bfP(W_{n\sigma^2_{X}}>y)\leq c_s \exp\{\frac{-c_s'y^2}{n\sigma^2_{X}}\}, \forall y>0.$
Hence for any $\varepsilon>0$ we can find constants $c_s'', \epsilon>0$ such that
\beqlb\label{lem4.1a}
\bfP\left(\left|\sum_{i=1}^{n}\bar{X}_i\right|\geq 2\varepsilon n\right)\leq\bfP\left(\left|\sum_{i=1}^{n}\bar{X}_i\right|\geq 2\varepsilon n, |W_{n\sigma^2_{X}}|\leq \varepsilon n\right)+\bfP\left(|W_{n\sigma^2_{X}}|> \varepsilon n\right)
\leq c_s''n^{-1-\epsilon}.\no
\eeqlb
Note that $\bfP(|n^{-1}\sum_{i=t_n+1}^{t_n+n}X_i-\bfE(X_1)|\geq 2\varepsilon)=\bfP\left(\left|\sum_{i=1}^{n}\bar{X}_i\right|\geq 2\varepsilon n\right).$ Hence we get \eqref{l4l} by Borel-Cantelli 0-1 law and the above inequalities.

\noindent{\bf Proof of Corollary \ref{mog0}}

With the help of \eqref{mogc1u}, to prove Corollary \ref{mog0} we only need to show
\begin{eqnarray}\label{mog.1>}\varliminf \limits_{n\rightarrow +\infty}\frac{\log \mbfP_\L
\left(\substack{\forall_{0\leq i\leq n}~\frac{T_{t_n+i}}{n^\alpha}\in \left[g\left(\frac{i}{n}\right),h\left(\frac{i}{n}\right)\right],~\\ \frac{T_{t_n+n}}{n^\alpha}\in [a',b'],~~~\xi_{t_n+i}\leq e^{n^{v}}}\Bigg|T_{t_n}=h(0)n^{\alpha}\right)}{n^{1-2\alpha}}
\geq -C_{g,h}\gamma_{\sigma},~~{\rm \mathbf{P}-a.s.}\end{eqnarray}

First we should note that Condition 4 means that
$\mbfE\left(|\log\mbfP_\L(T_{1}\leq y|T_{0}=0)|^{2+\epsilon}\right)<+\infty$ by the many-to-one formula. Therefore, according to the Monotone convergence theorem,
there exist constants $x<0, A>0$ such that
\begin{eqnarray}\label{mog.cc}\mbfE\left(|\log\mbfP_\L(x\leq T_{1}\leq y, \xi_1\leq A|T_{0}=0)|^{2+\epsilon}\right)<+\infty.\end{eqnarray}
Choose a small enough constant~$\delta>0$ such that~$\frac{g(0)-h(0)}{x}>\delta,$ then there exists a constant~$\epsilon>0$ such that~$g(0)+\epsilon-h(0)< x\delta.$ Denote~$\delta_n:=\lfloor\delta n^{\alpha}\rfloor.$ By the continuity of~$g$ 
one can choose $n$ large enough such that~$g(s)\leq g(0)+\epsilon$ for any $s\in[0,\frac{\delta_n}{n}].$
~Note that~$h(s)>h(0), s\in(0,1],$ then for any~$i\in[0,\delta_n]\cap\bfN,$ we have
\begin{eqnarray}\label{gh}
(g(i/n)-h(0))n^{\alpha}< xi\leq yi\leq (h(i/n)-h(0))n^{\alpha}\end{eqnarray} and $(h(\delta_n/n)-h(0))n^{\alpha}>0.$
\eqref{gh} can be derived from
~$$(g(i/n)-h(0))n^{\alpha}\leq (g(0)+\epsilon-h(0))n^{\alpha}< x\delta n^{\alpha}\leq xi, ~\forall i\in[0,\delta_n]\cap\bfN.$$ By Markov property, we have
\begin{eqnarray}\label{mog.2}
 &&\mbfP_\L
\left(\substack{\forall_{0\leq i\leq n}~\frac{T_{t_n+i}}{n^\alpha}\in \left[g\left(\frac{i}{n}\right),h\left(\frac{i}{n}\right)\right],~\\ \frac{T_{t_n+n}}{n^\alpha}\in [a',b'],~~~\xi_{t_n+i}\leq e^{n^{v}}}\Bigg|T_{t_n}=h(0)n^{\alpha}\right)\nonumber
\\&\geq&\mbfP_\L
(\forall_{0\leq i\leq \delta_n}~T_{t_n+i}\in [xi,yi],~\xi_{t_n+i}\leq e^{n^{v}}|T_{t_n}=0)\nonumber
\\&\times&\inf_{z\in[x\delta_n,y\delta_n]}\mbfP_\L
\left(\substack{\forall_{\delta_n\leq i\leq n}~\frac{T_{t_n+i}}{n^\alpha}\in [g(\frac{i}{n})-h(0),h(\frac{i}{n})-h(0)],\\ \frac{T_{t_n+n}}{n^\alpha}\in [a'-h(0),b'-h(0)],~\xi_{t_n+i}\leq e^{n^{v}}}\Bigg|T_{t_n+\delta_n}=z\right)\nonumber
\\&\geq&\prod_{m=1}^{\delta_n}\mbfP_\L
(T_{t_n+m}\in [x,y],\xi_{t_n+m}\leq e^{n^{v}}|T_{t_n+m-1}=0)\nonumber
\\&\times&\inf_{z\in[x\delta_n,y\delta_n]}\mbfP_\L
\Big(\substack{\forall_{0\leq i\leq n-\delta_n}~\frac{T_{t_n+\delta_n+i}}{n^\alpha}\in [g(\frac{i+\delta_n}{n})-h(0),h(\frac{i+\delta_n}{n})-h(0)],\\ \frac{T_{t_n+n}}{n^\alpha}\in [a'-h(0),b'-h(0)],~\xi_{t_n+\delta_n+i}\leq e^{n^{v}}}\Big|T_{_{t_n+\delta_n}}=z\Big).
\end{eqnarray}
We observe that for any $i\leq n-\delta_n$,
\begin{eqnarray}\label{uniform}
\left|\frac{i}{n-\delta_n}-\frac{i+\delta_n}{n}\right|=\left|\frac{\delta_n^2+\delta_ni-n\delta_n}{(n-\delta_n)n}\right|\leq\frac{\delta_n^2+n\delta_n}{(n-\delta_n)n}\leq\frac{\delta_n^2+\delta_n}{n-\delta_n}.
\end{eqnarray}
By recalling $\delta_n:=\lfloor\delta n^{\alpha}\rfloor$ and $\alpha\in(0,\frac{1}{3}],$ one sees that $\frac{\delta_n^2+\delta_n}{n-\delta_n}\ra 0.$
Moreover, note that $g,h$ are both continuous functions on~$[0,1]$ hence they are bounded and uniformly continuous. Thus for any given~$\varepsilon>0$ and $n$ large enough, from \eqref{uniform} we have $$\sup_{i\leq n-\delta_n}\left|\frac{(n-\delta_n)^\alpha}{n^\alpha}g\left(\frac{i}{n-\delta_n}\right)-g\left(\frac{i+\delta_n}{n}\right)\right|<\frac{\varepsilon}{4},~~~\frac{n^\alpha}{(n-\delta_n)^\alpha}<2,$$
$$ \sup_{i\leq n-\delta_n}\left|\frac{(n-\delta_n)^\alpha}{n^\alpha}h\left(\frac{i}{n-\delta_n}\right)-h\left(\frac{i+\delta_n}{n}\right)\right|<\frac{\varepsilon}{4} ~\text{and}~ \left|\frac{(n-\delta_n)^\alpha}{n^\alpha}h(0)-h(0)\right|<\frac{\varepsilon}{4}.$$
From the above analysis, we can obtain
\begin{eqnarray}\label{mog.3}
&&\frac{\log\inf\limits_{z\in[x\delta_n,y\delta_n]}\mbfP_\L
\Big(\substack{\forall_{i\leq n-\delta_n}~\frac{T_{t_n+\delta_n+i}}{n^\alpha}\in [g(\frac{i+\delta_n}{n})-h(0),h(\frac{i+\delta_n}{n})-h(0)],\\ \frac{T_{t_n+n}}{n^\alpha}\in [a'-h(0),b'-h(0)],~\xi_{t_n+\delta_n+i}\leq e^{n^{v}}}\Big|T_{_{t_n+\delta_n}}=z\Big)}{n^{1-2\alpha}}\nonumber
\\&\geq&\frac{\log\inf\limits_{z\in[x\delta_n,y\delta_n]}\mbfP_\L
\Big(\substack{\forall_{i\leq n-\delta_n}~\frac{T_{t_n+\delta_n+i}}{(n-\delta_n)^\alpha}+h(0)\in [g(\frac{i}{n-\delta_n})+\varepsilon,h(\frac{i}{n-\delta_n})-\varepsilon],\\ \frac{T_{t_n+n}}{(n-\delta_n)^\alpha}\in [a'-h(0)+\varepsilon,b'-h(0)-\varepsilon],~\xi_{t_n+\delta_n+i}\leq e^{n^{v}}}\Big|T_{_{t_n+\delta_n}=z}\Big)\left(\frac{n-\delta_n}{n}\right)^{^{1-2\alpha}}}{(n-\delta_n)^{1-2\alpha}}\nonumber
\\&\geq& -C_{g+\varepsilon,h-\varepsilon}\gamma_{\sigma}.
\end{eqnarray}
The last inequality holds because one can view~$n-M$ (resp., $t_n+\delta_n$) as~the $n$ (resp., $t_n$) in Corollary \ref{mogc2} and hence we can apply Corollary \ref{mogc2} to get the lower bound $-C_{g+\varepsilon,h-\varepsilon}\gamma_{\sigma}$.

Note that for $n$ large enough we have
\begin{eqnarray}\label{mog.4}&& \varliminf_{n\ra\iy}\frac{\log\prod_{m=1}^{\delta_n}\mbfP_\L(T_{t_n+m}\in [x,y],\xi_{t_n+m}\leq e^{n^v}|T_{t_n+m-1}=0)}{n^{1-2\alpha}}\nonumber
\\&\geq& \varliminf_{n\ra\iy}\frac{\sum_{m=1}^{\delta_n}\log\mbfP_\L(T_{t_n+m}\in [x,y],\xi_{t_n+m}\leq A|T_{t_n+m-1}=0)}{\delta_n}\frac{\delta_n}{n^{1-2\alpha}}\nonumber
\\&=& \delta~\mbfE(\log\mbfP_\L(T_{1}\in [x,y],\xi_{1}\leq A|T_{0}=0)).\end{eqnarray}
The last equality is because of Lemma 4.1 and the fact that the random sequence
$$\{\mbfP_\L
(T_{t_n+m}\in [x,y],\xi_{t_n+m}\leq A|T_{t_n+m-1}=0)\}_{m\in\{1,2,...M\}}$$ is i.i.d. under~$\mathbf{P}.$
Note that we can use the law of large number but not Lemma \ref{lem4} when $t_n\equiv k$. That is why we only need $\lambda_5\geq 1$ when $t_n$ does not depend on $n$.
Recall that $\alpha=\frac{1}{3}$ hence~$\lim_{n\ra\iy}\frac{\delta_n}{n^{1-2\alpha}}=\delta.$
(Note that $\lim_{n\ra\iy}\frac{\delta_n}{n^{1-2\alpha}}\leq \delta$ as long as $\alpha\in(0,1/3].$ That is why we say Corollary \ref{mog0} still holds as long as $\alpha\in(0, 1/3]$.) 
Letting~$\varepsilon\ra 0, \delta\ra 0,$  we get Corollary~\ref{mog0} according to \eqref{mog.2}-\eqref{mog.4}. \qed

\bigskip
\section{Proofs of Propositions and Example}

\noindent{\bf Proof of Proposition \ref{propcon1}}

Recall the notation $\Phi_1=e^{\kappa_1}.$ Note that for any constant $a, A_i\in\bfR, n\in\bfN,$ by Cauchy-Schwarz inequality we have
\begin{eqnarray}\label{pfpro1}
\left(\sum^n_{i=1}A^2_ie^{aA_i}\right)\left(\sum^n_{i=1}e^{aA_i}\right)
&=&\left(\sum^n_{i=1}A^2_ie^{2aA_i}+\sum_{1\leq i<j\leq n}(A^2_i+A^2_j)e^{a(A_i+A_j)}\right)\no
\\&\geq&\left(\sum^n_{i=1}A^2_ie^{2aA_i}+\sum_{1\leq i<j\leq n}2A_iA_je^{a(A_i+A_j)}\right)\no
\\&=&\left(\sum^n_{i=1}A_ie^{aA_i}\right)^2.
\end{eqnarray}
By the H\"{o}lder inequality and \eqref{pfpro1} we get
\begin{eqnarray}\label{pfpro2}
\Phi''_1(\vartheta)\Phi_1(\vartheta)
&\geq&\left[\mbfE_{\L}\left(\sqrt{\left(\sum_{i=1}^{N(\phi)}\zeta_i^2(\phi)e^{-\vartheta\zeta_{i}(\phi)}\right)\left(\sum_{i=1}^{N(\phi)}e^{-\vartheta\zeta_{i}(\phi)}\right)}\right)\right]^2\no
\\&\geq&\left[\mbfE_{\L}\left(\left|\sum_{i=1}^{N(\phi)}\zeta_i(\phi)e^{-\vartheta\zeta_{i}(\phi)}\right|\right)\right]^2\no
\\&\geq&[\Phi'_1(\vartheta)]^2.
\end{eqnarray}
Note that
$$\kappa_1''(\vartheta)=\frac{\mbfE_{\L}\left(\sum_{i=1}^{N(\phi)}\zeta^2_i(\phi)e^{-\vartheta \zeta_i(\phi)}\right)\mbfE_{\L}\left(\sum_{i=1}^{N(\phi)}e^{-\vartheta \zeta_i(\phi)}\right)-\left[\mbfE_{\L}\left(\sum_{i=1}^{N(\phi)}\zeta_i(\phi)e^{-\vartheta \zeta_i(\phi)}\right)\right]^2}{\left[\mbfE_{\L}\left(\sum_{i=1}^{N(\phi)}e^{-\vartheta \zeta_i(\phi)}\right)\right]^2}.$$
Hence by \eqref{pfpro2} we know $\kappa''_1(\vartheta)\geq0, \mathbf{P}-{\rm a.s.}$

We call the environment $L$ a constant jumping environment if we can find a constant $c(L)$ such that $\mbfP_\L(\zeta_{i}(\phi)=c(L),\forall i\leq N(\phi)|\L=L)=1,$ where we write $c(L)$ but not $c$ since different environments $L$ may correspond to different constants.

Now we use the proof by contradiction to show that $\kappa''_1(\vartheta)|_{\L=L}>0$ when $L$ is not a constant jumping environment.

We suppose that $\kappa''_1(\vartheta)=0,$ which means that the three $``\geq"$ in \eqref{pfpro2} are all $``=".$ By the H{\"o}lder's inequality, we
know that the first $``\geq"$ can be $``="$ only if
 there exists a constant $b\geq0$ ($b$ may depend on the realization of $\L$) such that
 \begin{eqnarray}\label{pfpro3}
 \left(\sum_{i=1}^{N(\phi)}\zeta_i^2(\phi)e^{-\theta\zeta_{i}(\phi)}\right)=b\left(\sum_{i=1}^{N(\phi)}e^{-\theta\zeta_{i}(\phi)}\right),~~{\rm \mbfP_{\L}-a.s.} \end{eqnarray}
 According to \eqref{pfpro1}, the second inequality in \eqref{pfpro2} holds only if for any $n\geq 1,$ on the event $\{N(\phi)=n\}$ we have \begin{eqnarray}\label{pfpro4}\zeta_{1}(\phi)=\zeta_{2}(\phi)=\cdots=\zeta_{n}(\phi)(:=b_n),~\text{where}~ b_n~ \text{is a constant}.\end{eqnarray}
 Combining \eqref{pfpro3} with \eqref{pfpro4} we deduce that $$b_n=\sqrt{b}~ \text{or}~-\sqrt{b}, ~\forall n\geq 1.$$  Note that for any random variable $Y,$ $|\bfE Y|=\bfE|Y|$ only if $Y$ is non-negative or non-positive. Therefore, from the third inequality in \eqref{pfpro2} we deduce that
 $$b_n=\sqrt{b}, ~\forall n\geq 1~~\text{or}~~ b_n=-\sqrt{b}, ~\forall n\geq 1,$$
which means that we can find a constant $c(L)$ ($-\sqrt{b}$ or $\sqrt{b}$) such that
$$\mbfP_{\L=L}(\zeta_{i}(\phi)=c(L),\forall i\leq N(\phi))=1.$$ This contradicts the assumption that $L$ is not a constant jumping environment. So far we have shown that
$\kappa''_1(\vartheta)|_{\L=L}>0$ when $L$ is not a constant jumping environment.

Next, we use proof by contradiction again to show $\mathbf{P}(\kappa''_1(\vartheta)>0)>0.$ We suppose that $\mathbf{P}(\kappa''_1(\vartheta)>0)=0,$ which is equivalent to say $\mathbf{P}(\kappa''_1(\vartheta)=0)=1.$~According to the above conclusion, it means that in the sense of ${\rm \mbfP-a.s.},$  for any realization $L$ of $\L$ we can find a constant $c(L)$ such that $\mbfP_{\L=L}(\zeta_{i}(\phi)=c(L),\forall i\leq N(\phi))=1.$ In this case, it is easy to see $$\kappa_1(\vartheta)=-\vartheta c(\L)+\log\mbfE_{\L}\left(N(\phi)\right)~ ~\text{and}~~\kappa'_1(\vartheta)=\frac{-c(\L)e^{-\vartheta c(\L)}\mbfE_{\L}\left(N(\phi)\right)}{e^{-\vartheta c(\L)}\mbfE_{\L}\left(N(\phi)\right)}=-c(\L).$$ Recall that in Condition 1 we have assumed that $\kappa(\vartheta)-\vartheta\kappa'(\vartheta)=0.$ Hence we have $$0=\kappa(\vartheta)-\vartheta\kappa'(\vartheta)=\mbfE(\kappa_1(\vartheta)-\vartheta\kappa'_1(\vartheta))=\mbfE(\log\mbfE_{\L}\left(N(\phi)\right))=\kappa(0).$$
But this contradicts the assumption $\kappa(0)>0$ in Condition 1. Hence we get $\mathbf{P}(\kappa''_1(\vartheta)>0)>0.$~\qed



{\noindent\bf Proof of Proposition \ref{propcon2}}

First we show \eqref{c21a} $\Rightarrow$ \eqref{c2'b}.

By the many-to-one formula, \eqref{c2'b} is equal to $\mbfE((\mbfE_{\L}|T_1-\mbfE_{\L}T_1|^{\lambda_2})^{\lambda_1})<+\infty.$ By Jensen's inequality we see that if there exists $\lambda_6>\frac{3}{2}$ such that $\mbfE((\mbfE_{\L}(T_1-\mbfE_{\L}T_1)^{4})^{\lambda_6})<+\infty,$ then $\mbfE((\mbfE_{\L}(T_1-\mbfE_{\L}T_1)^{\lambda_2})^{\lambda_1})<+\infty$ since $\lambda_2>2, \lambda_1>3.$ Moreover, using the many-to-one formula again we get
$$\mbfE_{\L}((T_1-\mbfE_{\L}T_1)^{4})=\frac{\mbfE_{\L}\left(\sum_{i=1}^{N(\phi)}|\zeta_i(\phi)+\kappa'_1(\vartheta)|^{4}e^{-\vartheta \zeta_i(\phi)}\right)}{\mbfE_{\L}\left(\sum_{i=1}^{N(\phi)}e^{-\vartheta \zeta_i(\phi)}\right)}.$$
Now we only need to show
\begin{eqnarray}\label{pro2.2}\frac{\mbfE_{\L}\left(\sum_{i=1}^{N(\phi)}|\zeta_i(\phi)+\kappa'_1(\vartheta)|^{4}e^{-\vartheta \zeta_i(\phi)}\right)}{\mbfE_{\L}\left(\sum_{i=1}^{N(\phi)}e^{-\vartheta \zeta_i(\phi)}\right)}=\kappa^{(4)}_1(\vartheta)+3[\kappa''_1(\vartheta)]^2,\end{eqnarray}
which will completes this proof.

Let $\Phi_1:=e^{\kappa_1}.$ Obviously, $\mbfE_{\L}\left(\sum_{i=1}^{N(\phi)}\zeta^n_i(\phi)e^{-\vartheta \zeta_i(\phi)}\right)=(-1)^n\Phi^{(n)}_1.$ Moreover, we can see
$$\Phi'_1=\kappa'_1\Phi_1,\Phi''_1=([\kappa'_1]^2+\kappa''_1)\Phi_1,~\Phi^{(3)}_1=((\kappa'_1)^3+3\kappa'_1\kappa''_1+\kappa^{(3)}_1)\Phi_1.$$
$$\Phi^{(4)}_1=((\kappa'_1)^4+6(\kappa'_1)^2\kappa''_1+4\kappa'_1\kappa^{(3)}_1+3(\kappa''_1)^2+\kappa^{(4)}_1)\Phi_1.$$

Then we get \eqref{pro2.2} by direct calculation.

Secondly, we show \eqref{c22a} $\Rightarrow$ \eqref{c21a}.

Without loss of generality we take $\lambda_6=2$ and hence we only need to show that $\mbfE(S_1^{8}+(\kappa'_1)^8)<+\infty.$

Note that there exists a $c_1$ such that $|x|\leq c_1e^{\frac{\lambda_7}{8}x}+c_1e^{-\frac{\lambda_8}{8}x}, \forall x\in\bfR.$ Therefore,
$$|\kappa'_1|\leq \frac{\mbfE_{\L}\left(\sum_{i=1}^{N(\phi)}|\zeta_i(\phi)|e^{-\vartheta \zeta_i(\phi)}\right)}{\mbfE_{\L}\left(\sum_{i=1}^{N(\phi)}e^{-\vartheta \zeta_i(\phi)}\right)}\leq c_1e^{\kappa_1(\theta-\frac{\lambda_7}{8})-\kappa_1(\vartheta)}+c_1e^{\kappa_1(\theta+\frac{\lambda_8}{8})-\kappa_1(\vartheta)}.$$
By the convexity of $\kappa_1$ we see $\kappa_1(\theta-\frac{\lambda_7}{8})-\kappa_1(\vartheta)\leq\kappa_1(\theta-(i+1)\frac{\lambda_7}{8})-\kappa_1(\vartheta-i\frac{\lambda_7}{8})$ and~
$\kappa_1(\theta+\frac{\lambda_8}{8})-\kappa_1(\vartheta)\leq\kappa_1(\theta+(i+1)\frac{\lambda_7}{8})-\kappa_1(\vartheta+i\frac{\lambda_7}{8})$.
Note that $(a+b)^8\leq 2^7(a^8+b^8), \forall a,b\in\bfR.$ Hence  $$\mbfE((\kappa'_1)^8)\leq c_12^7\mbfE(e^{\kappa_1(\theta-\lambda_7)-\kappa_1(\vartheta)}+e^{\kappa_1(\theta+\lambda_8)-\kappa_1(\vartheta)})<+\infty.$$
The way to show $\mbfE(S_1^8)<+\infty$ is similar. We can also find a $c_2>0$ such that $x^8\leq c_2(e^{\lambda_7x}+e^{-\lambda_8x}), \forall x\in\bfR.$ By the many-to-one formula we see
$$\mbfE_{\L}(e^{\lambda_7S_1}+e^{-\lambda_8S_1})=e^{\kappa_1(\theta-\lambda_7)-\kappa_1(\vartheta)}+e^{\kappa_1(\theta+\lambda_8)-\kappa_1(\vartheta)}.$$
So far we have shown $\mbfE(S_1^{8}+(\kappa'_1)^8)<+\infty.$ \qed

\noindent{\bf Proof of Proposition \ref{propcon3}}

It is obvious by Jensen's inequality and Proposition \ref{propcon2}.

\noindent{\bf Proof of Remark \ref{5--}}

Condition 3 is set in order to ensure that Corollary \ref{mogc01} holds (and thus Corollary \ref{mogc2} holds). Here we will show if the branching is independent of the displacement, then Corollary \ref{mogc01} still holds even though the Condition 3 only contains $\mbfE([\log^+\mbfE_{\L}(N(\phi)^{1+\lambda_4})]^{\lambda_3})<+\infty.$  Checking the proof of Corollary \ref{mogc01}, we can obtain the upper bound of $\mbfP_{\L}(\xi_i>e^{n^v})$ according to the independence of $N(u)$ and $\zeta_i(u)$ but not according to the H\"{o}lder's inequality like what we do in \eqref{Holder}. By the independence we have $$\mbfP_{\L}(\xi_i>e^{n^v})=\frac{\mbfE_{\L}\left(\1_{\{N(u)>e^{n^v}\}}N(u)\right)}{\mbfE_{\L}\left(N(u)\right)}\leq  \frac{\mbfE_{\L}(N(u)^{1+\lambda_4})}{e^{\lambda_4n^v}\mbfE_{\L}\left(N(u)\right)}, ~~|u|=i-1.$$
Hence we can adjust the definition of $\hat{I}_n$ (see \eqref{tilIn}) from $$\hat{I}_n:=\left\{\max_{|u|\leq n}\mbfE_{\L}(N(u)^{1+\lambda_3})\leq e^{n^{v_2}}, \max_{i\leq n}[(1-v_1)\kappa_i(\vartheta+\lambda_4)-\kappa_i(\vartheta)]\leq n^{v_2}\right\}$$ to $\hat{I}_n:=\left\{\max_{|u|\leq n}\mbfE_{\L}(N(u)^{1+\lambda_3})\leq e^{n^{v_2}}\right\}.$ Therefore, from \eqref{PROP3} to the end of the proof of Corollary \ref{mogc01} (i.e. to get the convergence in Corollary \ref{mogc01}), we only need to dominate the tail of $\mbfE_{\L}(N(u)^{1+\lambda_3})$. \qed

\noindent{\bf Proof of~Example \ref{exb0}}

In this proof, $N(\mathfrak{m}), \mu(\mathfrak{m})$ and $\sigma^2(\mathfrak{m})$ are always abbreviated by $N, \mu$ and $\sigma$ respectively. Hence we should note that $N, \mu, \sigma$ are random variables under the annealed law $\mbfP.$
By the statement of this example we see
$\kappa_1(x)=\log\mbfE_{\L}N-x\mu+\frac{1}{2}x^2\sigma^2,\forall x\in\bfR$ hence
$\kappa'_1(x)=\sigma^2 x-\mu,~\kappa''_1(x)=\sigma^2.$ Note that \eqref{exb1} implies $(\kappa(0)=)\mbfE(\log\mbfE_{\L}N)\in(0,\infty)$ and \eqref{exb2} implies $\mbfE(\sigma^2)>0.$
Taking $\vartheta:=\sqrt{\frac{2\mbfE(\log\mbfE_{\L}N)}{\mbfE(\sigma^2)}},$ we have $\vartheta\kappa'(\vartheta)-\kappa(\vartheta)=0,$ which means that Condition 1 holds.

Obviously, \eqref{exb1} and \eqref{exb2} imply that \eqref{c2'a} holds. Note that $\kappa^{(4)}_1(x)\equiv0,$ hence \eqref{c2'b} also holds by Proposition 2.2 and \eqref{exb2}.  Thus the example satisfies Condition 2.

Condition 3 holds because of Remark \ref{5--} and \eqref{exb1}.

At last, we verify the Condition 4. By the proof of Corollary \ref{mog0} we know that Condition 4 is set to ensure $\mbfE\left(|\log\mbfP_\L(T_{1}\leq y|T_{0}=0)|^{\lambda_5}\right)<+\infty$ for some $\lambda_5>2$. From the continuity of normal distribution and the arbitrary of $y$ in Condition 4 we can see it is enough to show \begin{eqnarray}\label{mog.cex}\mbfE\left(|\log\mbfP_\L(T_{1}\leq 0|T_{0}=0)|^{\lambda_5}\right)<+\infty.\end{eqnarray} Next we devote to prove \eqref{mog.cex}. By many-to-one formula we see that for any $\lambda\in\bfR,$
\begin{eqnarray}
\log\mbfE_{\L}(e^{\lambda(T_1-\mbfE_{\L}T_1)})&=&\kappa_1(\vartheta-\lambda\vartheta)-\kappa_1(\vartheta)+\lambda\vartheta\kappa'_1(\vartheta)
=\frac{1}{2}\lambda^2\vartheta^2\sigma^2,~{\rm \mathbf{P}-a.s.}\nonumber
\end{eqnarray}
That is to say, under $\mbfP_{\L},$ $T_1-\mbfE_{\L}T_1$ has the normal distribution $\mathcal{N}(0,\vartheta^2\sigma^2).$ Hence we have \begin{eqnarray}\label{ex+1}\mbfP_{\L}(T_1\leq 0)&=&\mbfP_{\L}(T_1-\mbfE_{\L}T_1\leq -\mbfE_{\L}T_1)\no\\&\geq&\frac{1}{\sqrt{2\pi}\vartheta\sigma}e^{-\frac{(\mbfE_{\L}T_1+1)^2}{2\vartheta^2\sigma^2}}\1_{\{\mbfE_{\L}T_1\geq 0\}}+\frac{1}{2}\1_{\{\mbfE_{\L}T_1< 0\}}.
\end{eqnarray}
Note that $\mbfE_{\L}T_1=\kappa_1(\vartheta)-\vartheta\kappa'_1(\vartheta)=\log\mbfE_{\L}N-\frac{1}{2}\vartheta^2\sigma^2$ thus $\exp \{-\frac{(\mbfE_{\L}T_1+1)^2}{2\vartheta^2\sigma^2}\}\1_{\{\mbfE_{\L}T_1\geq 0\}}\geq \exp \{-\frac{(\log\mbfE_{\L}N+1)^2}{2\vartheta^2\sigma^2}\}\1_{\{\mbfE_{\L}T_1\geq 0\}}.$ 
Then we have
\begin{eqnarray}\label{ex+2}
&&\mbfE\left(|\log\mbfP_\L(T_{1}\leq 0|T_{0}=0)|^{\lambda_5}\1_{\{\mbfE_{\L}T_1\geq 0, \sigma\leq \frac{1}{\sqrt{2\pi}\vartheta}\}}\right)\leq \mbfE\left(\frac{(\log\mbfE_{\L}N+1)^{2\lambda_5}}{(2\vartheta^2\sigma^2)^{\lambda_5}}\right).
\end{eqnarray}
It is known that we can find a constant $c>0$ such that $cx^2\geq \log(\sqrt{2\pi}\vartheta x)$ when $\sqrt{2\pi}\vartheta x>1.$ Note that on $\{\mbfE_{\L}T_1\geq 0\},$ we have $\log(\sqrt{2\pi}\vartheta\sigma)+\frac{(\mbfE_{\L}T_1+1)^2}{2\vartheta^2\sigma^2}\geq 0$ and $\log\mbfE_{\L}N=\log^+\mbfE_{\L}N$. Hence on the event $\{\mbfE_{\L}T_1\geq 0\}\cap\{\sqrt{2\pi}\vartheta \sigma>1\},$ we have
\begin{eqnarray}0&\leq& \log(\sqrt{2\pi}\vartheta\sigma)+\frac{(\mbfE_{\L}T_1+1)^2}{2\vartheta^2\sigma^2}\leq c\sigma^2+\frac{(\log\mbfE_{\L}N_1+1)^2}{2\vartheta^2\sigma^2}\no
\\&\leq&\frac{2c\log \mbfE_{\L}N}{\vartheta^2}+\frac{(\log\mbfE_{\L}N_1+1)^2}{2\vartheta^2\sigma^2}\leq \frac{4c(\log \mbfE_{\L}N)^2}{\vartheta^4\sigma^2}+\frac{(\log\mbfE_{\L}N+1)^2}{2\vartheta^2\sigma^2}.\no\end{eqnarray}
From the above analysis, we can find a constant $c_1>0$ such that
\begin{eqnarray}\label{ex+3}
\mbfE\left(|\log\mbfP_\L(T_{1}\leq 0|T_{0}=0)|^{\lambda_5}\1_{\{\mbfE_{\L}T_1\geq 0, \sigma> \frac{1}{\sqrt{2\pi}\vartheta}\}}\right)&\leq& c_1\mbfE\left(\frac{(\log\mbfE_{\L}N+1)^{2\lambda_5}}{(2\vartheta^2\sigma^2)^{\lambda_5}}\1_{\{\mbfE_{\L}T_1\geq 0, \sigma> \frac{1}{\sqrt{2\pi}\vartheta}\}}\right)\no
\\&\leq& c_1\mbfE\left(\frac{(\log\mbfE_{\L}N+1)^{2\lambda_5}}{(2\vartheta^2\sigma^2)^{\lambda_5}}\right).
\end{eqnarray}
Combining with \eqref{ex+1}-\eqref{ex+3} and the independence of $N$ and $\sigma$ we have
\begin{eqnarray}
&&\mbfE\left(|\log\mbfP_\L(T_{1}\leq 0|T_{0}=0)|^{\lambda_5}\right)\leq(1+c_1)\mbfE\left((\log\mbfE_{\L}N+1)^{2\lambda_5}\right)\mbfE\left(\frac{1}{(2\vartheta^2\sigma^2)^{\lambda_5}}\right)+\lambda_5\log 2.\no\end{eqnarray}
Choosing $\lambda_5\in(2,\min\{3,\frac{\tau_2}{2}\}),$ then we verify \eqref{mog.cex} (and thus Condition 4) by \eqref{exb1}, \eqref{exb2} and the above inequality.\qed

\section{Proof of Theorems}




Let us sort out the relationship between Theorem \ref{1} and Theorem \ref{2}. It is obvious that {\rm(2a)}, {\rm(2b)}, and {\rm(2d)} imply {\rm(1c)}, {\rm(1d)} and {\rm(1b)} respectively. Moreover, based on the small deviation principle for RWre (Corollaries \ref{mogc01}-\ref{mog0}), the proofs of {\rm(2a)}, {\rm(2b)} and {\rm(2c)} will be very similar to the proofs of \cite[Proposition 1.5]{BJ2012}, \cite[Theorem 1.2]{AJ2011} and
\cite[Lemma 4.6]{GHS2011} respectively. Hence we omit these three proofs and only give the proofs of {\rm(1a)}, {\rm(1b)} and {\rm(2d)}. The main task in this section is to show {\rm(1b)}. As long as we get {\rm(1b)}, {\rm(1a)} is not hard to show by combining {\rm(1b)} and {\rm(2c)}. In addition, during our proof of {\rm(1b)}, we can find we have even shown {\rm(2d)}. We borrow some ideas from the proof of
\cite[Proposition 1.4]{BJ2012} to prove {\rm(1b)}, but there are some differences in details. (The main differences appear in \eqref{bu1}-\eqref{bu1+} and \eqref{e4.3.31}. For the convenience of reading and a better understanding, we give the complete proof of {\rm(1b)}.) We point out that even though the following long proof has some similarities with the proof in \cite{BJ2012}, it is far from the complete proof of the main theorems. The preparation in Section 4 and Proposition \ref{propcon1}, which are the new difficulties as we deal with the random environment, are also necessary to prove the main theorems.



\noindent{\bf Proofs of Theorem \ref{1} (1b)}

In this proof we let $M\in \bfN^+~,k\in \bfN$ 
and the barrier function $$\varphi_{\L}(i):=-\vartheta^{-1}K_i+ai^{\frac{1}{3}}.$$
First we should emphasize again that the $\varphi_{\L}(i)$ is a random variable depending on the random environment $\L$ because the $\{K_n\}_{n\in\bfN}$ is a random walk under $\mathbf{P}.$ To simplify the presentation, we sometimes omit $``{\rm \mathbf{P}-a.s.}"$ after some equalities or inequalities without causing confusion.
As the model of BRWre with barrier, in the $i$-th generation, any individual born above the barrier $\varphi_{\L}(i)$ is removed and consequently does not reproduce. We only care about the surviving particles in this system, i.e., the particle $z$ satisfying $V(z_i)\leq \varphi_{\L}(i),\forall i\leq |z|.$ We pick a surviving individual $z$ in generation $M^k$ and consider
$U_k(z),$ which represents the number of surviving descendants of $z$ in generation $M^{k+1}$.
We see that under this barrier the rightmost position of the surviving particle in the $M^k$-th generation is no larger than $-\vartheta^{-1}K_{M^k}+aM^{\frac{k}{3}}$. Therefore, we can find a naturally lower bound of $U_k(z)$ by considering, instead of $z$, a virtual individual $\tilde{z}$ in the same generation $M^k$ but
positioned on the barrier at $V(\tilde{z}):=-\vartheta^{-1}K_{M^k}+aM^{\frac{k}{3}}\geq V(z)$. Since the number and displacements of the descendants of $\tilde{z}$ are exactly the same as those of $z$, the descendants of $\tilde{z}$ are more likely to cross the barrier and thus be killed, which means that
$U_k(z)\geq U_k(\tilde{z}).$ Let $r_k:=e^{k^v}, k\in\bfN^+,$ where $v\in(\frac{2}{\lambda_3},\frac{1}{3})$ and $\lambda_3$ is the one introduced in Condition 3. Further, we define
$$Z_{k,b}:=\sharp\Big\{u\in\mathcal{T}_{M^{(k+1)}},u>\tilde{z}: \forall M^{k}<i\leq M^{(k+1)},\substack{ V(u_i)\in[(a-b)i^{1/3}-\vartheta^{-1} K_{i},ai^{1/3}-\vartheta^{-1}K_{i}],\\ N(u_{i-1})\leq r_k}\Big\},$$
where the notations $\sharp, \mathcal{T}_{n}, N(\cdot), >$ have been defined in Section 1 and $b\in(0,a).$ The exact value of constant $b$ will be given later. 
It is obvious that $Z_{k,b}\leq U_k(\tilde{z}).$ Hence for any $k\in \bfN,$ we have $Z_{k,b}\leq U_k(\tilde{z})\leq U_k(z),~ {\rm \mathbf{P}-a.s.}$

Recall the definition $Y_n:=\sharp\{u\in\mathcal{T}_{n}, \forall i\leq n, V(u_i)\leq ai^{\frac{1}{3}}-\vartheta^{-1} K_{i}\}$ in Theorem \ref{2} and define  
\begin{eqnarray}\label{e4.3.5-}P_{n, \L}&:=&\mbfP_{\L}(\forall 1\leq k\leq n,Y_{M^k}\geq\eta\mbfE_{\L}(Z_{k-1,b})),\end{eqnarray}
where the constant $\eta\in(0,1).$
Then by Markov property and the definition of $P_{n, \L}$ we have
\begin{eqnarray}\frac{P_{n+1,\L}}{P_{n,\L}}&:=&\mbfP_{\L}(\forall 1\leq k\leq n+1,Y_{M^k}\geq\eta\mbfE_{\L}(Z_{k-1,b})|
\forall 1\leq k\leq n,Y_{M^k}\geq\eta\mbfE_{\L}(Z_{k-1,b}))\no
\\&\geq&1-\prod_{i=1}^{\lfloor\eta\mbfE_{\L}(Z_{n-1,b})\rfloor}\mbfP_{\L}(U_n(z^{(i)})< \eta\mbfE_{\L}(Z_{n,b}))\no
\\&\geq&1-\mbfP_{\L}(Z_{n,b}< \eta\mbfE_{\L}(Z_{n,b}))^{\lfloor\eta\mbfE_{\L}(Z_{n-1,b})\rfloor},\end{eqnarray}
where $\{z^{(i)}, i=1,2,\ldots,\lfloor\eta\mbfE_{\L}(Z_{n-1,b})\rfloor\}$ represents $\lfloor\eta\mbfE_{\L}(Z_{n-1,b})\rfloor$ different surviving individuals in the
$M^k$-th generation.
Denote~\begin{eqnarray}\label{AKL}A_{k,\L}:=\mbfP_{\L}(Z_{k,b}\geq \eta\mbfE_{\L}(Z_{k,b})),\end{eqnarray} then we have
\begin{eqnarray}\label{5.3e}P_{n,\L}&\geq& P_{1,\L}\prod_{k=1}^{n-1}(1-(1-A_{k,\L})^{\lfloor\eta\mbfE_{\L}(Z_{k-1,b})\rfloor})\no
\\&\geq& P_{1,\L}\prod_{k=1}^{n-1}(1-e^{-\lfloor\eta\mbfE_{\L}(Z_{k-1,b})\rfloor A_{k,\L}}).\end{eqnarray}
Moreover, it is not hard to see
\begin{eqnarray}\label{e4.3.5--}
~\mbfP_{\L}(\mathcal{S})\geq \lim\limits_{n\rightarrow +\infty}P_{n,\L}.
\end{eqnarray}
Note that for any $\eta\in(0,1),$ by the facts~$Y_M\geq Z_{0,b}$~and $Z_{0,b}\geq 0$ we can see
\begin{eqnarray}\label{e4.3.5++}P_{1,\L}=\mbfP_{\L}(Y_{M}\geq\eta\mbfE_{\L}(Z_{0,b}))\geq \mbfP_{\L}(Z_{0,b}\geq\eta\mbfE_{\L}(Z_{0,b}))>0,~~{\rm \mathbf{P}-a.s.}\end{eqnarray}
Therefore, if we can show
\begin{eqnarray}\label{e4.3.5}
\sum_{k=1}^{+\infty}e^{-\lfloor\eta\mbfE_{\L}(Z_{k-1,b})\rfloor A_{k,\L}}<+\infty, ~~{\rm \mathbf{P}-a.s.},
\end{eqnarray}
which will imply~$\prod_{k=1}^{+\infty}(1-e^{-\lfloor\eta\mbfE_{\L}(Z_{k-1,b})\rfloor A_{k,\L}})>0,~{\rm \mathbf{P}-a.s.}$,
then combining with~\eqref{5.3e}-\eqref{e4.3.5++} we see
\begin{eqnarray}\label{e4.3.5+}\mbfP_{\L}(\mathcal{S})\geq \lim\limits_{n\rightarrow +\infty}P_{n,\L}\geq P_{1,\L}\prod_{k=1}^{+\infty}(1-e^{-\lfloor \eta\mbfE_{\L}(Z_{k-1,b})\rfloor A_{k,\L}})>0,~{\rm \mathbf{P}-a.s.},\end{eqnarray}
which is exactly the conclusion in~Theorem 2.1~(1b).~

Hence in the rest part of this proof we only need to show \eqref{e4.3.5}.
That is to say~we want to find the lower bound of $A_{k,\L}$.
By the~Paley-Zygmund inequality, we have
\begin{eqnarray}\label{e4.3.5zy}
A_{k,\L}:=\mbfP_{\L}(Z_{k,b}\geq \eta\mbfE_{\L}(Z_{k,b}))\geq(1-\eta)^2\frac{[\mbfE_{\L}(Z_{k,b})]^2}{\mbfE_{\L}(Z^2_{k,b})},~~{\rm \mathbf{P}-a.s.}\end{eqnarray}
According to \eqref{e4.3.5zy}, first we try to get the upper bound of $\mbfE_{\L}(Z^2_{k,b})$ by using the second moment method. Recalling the important information we set for the particle $\tilde{z}:$ ~$|\tilde{z}|=M^k,~ V(\tilde{z})=aM^{\frac{k}{3}}-\vartheta^{-1}K_{M^{k}}.$
Define $$\Theta:=\left\{u\in \mathcal{T}:\substack{u>\tilde{z},~~~|u|\leq M^{k+1},~\forall M^k< i\leq |u|,~N(u_{i-1})\leq r_{k}~\\ ~V(u_i)\in[(a-b)i^{\frac{1}{3}}-\vartheta^{-1} K_{i},~ai^{\frac{1}{3}}-\vartheta^{-1} K_{i}]}\right\},$$
thus $Z_{k,b}$ also has the representation
$Z_{k,b}:=\sum_{|u|=M^{k+1}}\1_{\{u\in\Theta\}}.$
For any particle $v$, denote
\begin{eqnarray}\label{e4.3.6a}
Z^v_k(\Theta):=\1_{\{v\in\Theta\}}\left(\sum_{|u|=M^{k+1},u>v}\1_{\{u\in\Theta\}}\right).\end{eqnarray}
Let~$w$ be a child of~$v$, i.e.,~$w_{|w|-1}=v$ and set
$Z^v_k(\Theta,w):=\sum_{\left\{w':|w'|=|v|+1,w'>v,w'\neq w\right\}}Z^{w'}_k(\Theta),$
which stands for the number of the surviving descendants of $v$ in generation $M^{k+1}$ who are not the descendants of $w$.
Through the definitions above, for any particle $u$ in generation $M^{k+1}$, $Z_{k, b}$ can also be expressed as
$$Z_{k,b}:=\1_{\{u\in\Theta\}}+\sum\limits_{j=M^k}^{M^{k+1}-1}Z_k^{u_j}(\Theta,u_{j+1}).$$
On the other hand, note that~$Z^2_{k,b}=Z_{k,b}\left(\sum_{|u|=M^{k+1}}\1_{\{u\in\Theta\}}\right)=\sum_{|u|=M^{k+1}}Z_{k,b}\1_{\{u\in\Theta\}},$ hence we have
\begin{eqnarray}\label{e4.3.6}
Z^2_{k,b}-Z_{k,b}&=&\sum\limits_{|u|=M^{k+1}}\1_{\{u\in\Theta\}}(Z_{k,b}-\1_{\{u\in\Theta\}})\no
\\&=&\sum\limits_{|u|=M^{k+1}}\sum\limits_{j=M^k}^{M^{k+1}-1}\1_{\{u\in\Theta\}}Z_k^{u_j}(\Theta,u_{j+1})\no
\\&=&\sum\limits_{j=M^k}^{M^{k+1}-1}\sum\limits_{|u|=M^{k+1}}\1_{\{u\in\Theta\}}Z_k^{u_j}(\Theta,u_{j+1}).\end{eqnarray}
 We observe that for any two particles $u^{(1)},u^{(2)}$ in generation~$M^{k+1}$ with $u_j^{(1)}=u_j^{(2)}.$ If $u_{j+1}^{(1)}=u_{j+1}^{(2)}$ (i.e., they have a common ancestor in the $(j+1)$-th generation), 
then~$$Z_k^{u^{(1)}_j}(\Theta,u^{(1)}_{j+1})=Z_k^{u^{^{(1)}}_j}(\Theta,u^{(2)}_{j+1}).$$ Hence we have
\begin{eqnarray}\label{e4.3.7}\sum\limits_{|u|=M^{k+1}}\1_{\{u\in\Theta\}}Z_k^{u_j}(\Theta,u_{j+1})&=&\sum\limits_{|u'|=j+1}\sum\limits_{\substack{|u|=M^{k+1},\\u_{j+1}=u'}}\1_{\{u\in\Theta\}}Z_k^{u_j}(\Theta,u_{j+1})\no
\\&=&\sum\limits_{|u'|=j+1}\left(Z_k^{u'_j}(\Theta,u')\sum\limits_{\substack{|u|=M^{k+1},\\u_{j+1}=u'}}\1_{\{u\in\Theta\}}\right)\no
\\&=&\sum\limits_{|u'|=j+1}\left(Z_k^{u'_j}(\Theta,u')Z_k^{u'}(\Theta)\right).\end{eqnarray}
Combining \eqref{e4.3.7} with \eqref{e4.3.6} we obtain
\begin{eqnarray}\label{e4.3.8}
Z^2_{k,b}-Z_{k,b}&=&\sum\limits_{j=M^k}^{M^{k+1}-1}\sum\limits_{|u'|=j+1}\left(Z_k^{u'_j}(\Theta,u')Z_k^{u'}(\Theta)\right)
=\sum\limits_{|u'|=j+1}\sum\limits_{j=M^k}^{M^{k+1}-1}\left(Z_k^{u'_j}(\Theta,u')Z_k^{u'}(\Theta)\right)\no
\\&=&\sum\limits_{|v|=j}\sum\limits_{j=M^k+1}^{M^{k+1}}\left(Z_k^{\overleftarrow{v}}(\Theta,v)Z_k^{v}(\Theta)\right)=\sum\limits_{j=M^k+1}^{M^{k+1}}\sum\limits_{|v|=j}\left(Z_k^{\overleftarrow{v}}(\Theta,v)Z_k^{v}(\Theta)\right),
\end{eqnarray}
where $\overleftarrow{v}$ represents the parent of $v$ (i.e., $\overleftarrow{v}:=v_{|v|-1}$).

Now let us find the upper bound of $\sum\limits_{|v|=j}\left(Z_k^{\overleftarrow{v}}(\Theta,v)Z_k^{v}(\Theta)\right).$
Define the $\sigma-$algebra $\mathcal{F}_j:=\sigma(X(u),|u|<j).$ Then for any $j\in[M^k+1,M^{k+1}],$ we have
\begin{eqnarray}\label{e4.3.9}
\mbfE_{\L}\left(\sum\limits_{|v|=j}Z_k^{v}(\Theta)Z_k^{\overleftarrow{v}}(\Theta,v)\right)&=&\mbfE_{\L}\left(\mbfE_{\L}\Big(\sum\limits_{|v|=j}Z_k^{v}(\Theta)Z_k^{\overleftarrow{v}}(\Theta,v)|\mathcal{F}_{j}\Big)\right)\no
\\&=&\mbfE_{\L}\left(\sum\limits_{|v|=j}\mbfE_{\L}\Big[\sum_{v'=bro(v)}\left[Z_k^{v}(\Theta)Z_k^{v'}(\Theta)\right]\Big|\mathcal{F}_{j}\Big]\right),
\end{eqnarray}
where $~bro(v):=\{v':|v'|=|v|,v'_{|v|-1}=v_{|v|-1},v'\neq v\}$ and $v'\neq v$ represents that $v'$ is different from $v.$ That is to say, the set contains all siblings of $v.$ Note that $\sharp bro(v)$ is $\mathcal{F}_{j}$-measurable because of $|v|=j$. Therefore, we have
\begin{eqnarray}\label{e4.3.10}
\mbfE_{\L}\Bigg[\sum_{v'\in bro(v)}\left[Z_k^{v}(\Theta)Z_k^{v'}(\Theta)\right]\Big|\mathcal{F}_{j}\Bigg]&=&\sum_{v'\in bro(v)}\mbfE_{\L}\Big[Z_k^{v}(\Theta)Z_k^{v'}(\Theta)|\mathcal{F}_j\Big]\no
\\&=&\sum_{v'\in bro(v)}\mbfE_{\L}\Big[Z_k^{v}(\Theta)|\mathcal{F}_j\Big]\mbfE_{\L}\Big[Z_k^{v'}(\Theta)|\mathcal{F}_j\Big]\no
\\&=&\mbfE_{\L}\Big[Z_k^{v}(\Theta)|\mathcal{F}_j\Big]\sum_{v'\in bro(v)}\mbfE_{\L}\Big[Z_k^{v'}(\Theta)|\mathcal{F}_j\Big],
\end{eqnarray}
where the second equality is because conditionally on $\mathcal{F}_j$, $Z_k^{v}(\Theta)$ and~$Z_k^{v'}(\Theta)$ are independent of each other.
From~$\eqref{e4.3.9}$ and~$\eqref{e4.3.10}$ we see
\begin{eqnarray}\label{e4.3.11}
&&\mbfE_{\L}\Big(\sum\limits_{|v|=j}Z_k^{v}(\Theta)Z_k^{\overleftarrow{v}}(\Theta,v)\Big)\no
=\mbfE_{\L}\left\{\sum\limits_{|v|=j}\left( \mbfE_{\L}\Big[Z_k^{v}(\Theta)|\mathcal{F}_j\Big]\sum_{v'=bro(v)}\mbfE_{\L}\Big[Z_k^{v'}(\Theta)|\mathcal{F}_j\Big]\right)\right\}.
\end{eqnarray}
We see that if~$v\notin \Theta$, then~$\mbfE_{\L}\Big[Z_k^{v}(\Theta)\big|\mathcal{F}_j\Big]=0$ by~\eqref{e4.3.6a}. If~$v\in \Theta$, the set~$bro(v)$ has at most~$r_k$ ~elements because of the definition of $\Theta$. Moreover,~$\mbfE_{\L}\Big[Z_k^{v'}(\Theta)\big|\mathcal{F}_j\Big]$ only depends on~$v'$, hence we have
\begin{eqnarray}\label{e4.3.13}\mbfE_{\L}\Big(\sum\limits_{|v|=j}Z_k^{v}(\Theta)Z_k^{\overleftarrow{v}}(\Theta,v)\Big)
&\leq&\mbfE_{\L}\left\{\sum\limits_{|v|=j}\left( \mbfE_{\L}\Big[Z_k^{v}(\Theta)|\mathcal{F}_j\Big](r_k-1)\sup_{V(v')\in\mbr,|v'|=j}\mbfE_{\L}\Big[Z_k^{v'}(\Theta)\Big]\right)\right\}\no
\\&=&(r_k-1)\sup_{\substack{V(v')\in\mbr,|v'|=j}}\mbfE_{\L}\Big[Z_k^{v'}(\Theta)\Big]\mbfE_{\L}\left\{\sum\limits_{|v|=j}\left( \mbfE_{\L}\Big[Z_k^{v}(\Theta)|\mathcal{F}_j\Big]\right)\right\}\no
\\&=&(r_k-1)\sup_{V(v')\in\mbr,|v'|=j}\mbfE_{\L}\Big[Z_k^{v'}(\Theta)\Big]\mbfE_{\L}(Z_{k,b}).
\end{eqnarray}
The last equality is due to the smoothness of conditional expectation and the fact $Z_{k,b}=\sum\limits_{|v|=j}Z_k^{v}(\Theta).$ Taking expectations to both sides of $\eqref{e4.3.8}$ and then substituting it into $\eqref{e4.3.13},$ we get
$$\mbfE_{\L}(Z^2_{k,b})\leq\mbfE_{\L}(Z_{k,b})\left(1+(r_k-1)\sum_{j=M^k+1}^{M^{k+1}}\sup\limits_{|v'|=j, V(v')\in\mbr}\mbfE_{\L}\left[Z_k^{v'}(\Theta)\right]\right).$$
Combining with~\eqref{e4.3.5zy}, we get
\begin{eqnarray}\label{impc4}
A_{k,\L}\geq(1-\eta)^2\frac{[\mbfE_{\L}(Z_{k,b})]^2}{\mbfE_{\L}(Z^2_{k,b})}\no
\geq\frac{(1-\eta)^2\mbfE_{\L}(Z_{k,b})}{1+(r_k-1)\sum\limits_{j=M^k+1}^{M^{k+1}}\sup_{|v'|=j, V(v')\in\mbr}\mbfE_{\L}\Big[Z_k^{v'}(\Theta)\Big]},~~{\rm \mathbf{P}-a.s.}\end{eqnarray}
Next we begin to find the upper bound of~$\sup\limits_{|v'|=j,V(v')\in\mbr}\mbfE_{\L}\Big[Z_k^{v'}(\Theta)\Big]$.
By $\eqref{e4.3.6a}$ we know $$\sup_{\substack{|v'|=j, V(v')\in\mbr}}\Big[Z_k^{v'}(\Theta)\Big]=\sup_{\substack{|v'|=j, v'\in\Theta, V(v')\in\mbr}}\Big[Z_k^{v'}(\Theta)\Big].$$
According to the above definition~$I_j:=[(a-b)j^{\frac{1}{3}}-\vartheta^{-1} K_{j},~aj^{\frac{1}{3}}-\vartheta^{-1} K_{j}],$ we see that $v'\in\Theta$ means $V(v')\in I_j$ and hence
\begin{eqnarray}\label{e4.3.15}
&&\sup_{\substack{|v'|=j,~ V(v')\in\mbr}}\mbfE_{\L}\Big[Z_k^{v'}(\Theta)\Big]=\sup_{\substack{|v'|=j,~ V(v')\in I_j}}\mbfE_{\L}\left[\sum_{\substack{u_j=v',~|u|=M^{k+1}}}\1_{\{u\in\Theta\}}\right]\no
\\&\leq&\sup_{\substack{|v'|=j \\V(v')\in I_j}}\mbfE_{\L}\Bigg[\sum_{\substack{u_j=v'\\|u|=M^{k+1}}}\1_{\{\forall i\leq M^{k+1}-j,~V(u_{j+i})+\vartheta^{-1} K_{i+j}\in[(a-b)(i+j)^{\frac{1}{3}},~a(i+j)^{\frac{1}{3}}]\}}\Bigg].
\end{eqnarray}
By the definition of the shifted expectation $\mbfE^{j}_{\L}$ (see Section 3), we have
\begin{eqnarray}\label{e4.3.17}
&&\sup_{\substack{|v'|=j \\V(v')\in I_j}}\mbfE_{\L}\left[\sum_{\substack{u_j=v'\\|u|=M^{k+1}}}\1_{\{\forall i\leq M^{k+1}-j,~V(u_{j+i})\in[(a-b)(i+j)^{\frac{1}{3}}-\vartheta^{-1}K_{i+j},~a(i+j)^{\frac{1}{3}}-\vartheta^{-1}K_{i+j}]\}}\right]\no
\\&=&\sup_{y\in I_{j}}\mbfE^j_{\L}\left(\sum_{|v|=M^{k+1}-j}\1_{\{\forall i\leq M^{k+1}-j,~V(v_i)+\vartheta^{-1} K_{i+j}\in[-y+(a-b)(i+j)^{\frac{1}{3}},~-y+a(i+j)^{\frac{1}{3}}]\}}\right)\no
\\&=&\sup_{x\in [(a-b)j^{\frac{1}{3}},~aj^{\frac{1}{3}}]}\mbfE^j_{\L}\Big(e^{T_{M^{k+1}-j}}\1_{\{\forall i\leq M^{k+1}-j,~\frac{T_i}{\vartheta}\in[-x+(a-b)(i+j)^{\frac{1}{3}},~-x+a(i+j)^{\frac{1}{3}}]\}}\Big)\no
\\&\leq& e^{\vartheta aM^{\frac{k+1}{3}}-\vartheta(a-b)j^{\frac{1}{3}}}\times H_j,\end{eqnarray}
where $H_j:=\sup_{x\in [(a-b)j^{\frac{1}{3}},~aj^{\frac{1}{3}}]}\mbfP^j_{\L}\left(\forall_{i\leq M^{k+1}-j},
x+\frac{T_i}{\vartheta}\in[(a-b)(i+j)^{\frac{1}{3}},a(i+j)^{\frac{1}{3}}]\right).$
The second inequality above is due to the many-to-one formula \eqref{mto0} and the fact that $y+\vartheta^{-1} K_{j,\vartheta}\in[(a-b)j^{\frac{1}{3}},~aj^{\frac{1}{3}}].$

Now we divide the time axis $d_k:=M^{k+1}-M^{k}$ into $M^2-M$ segments equally. Let~$K(M):=M^2-M-1$ and $c_l:=M^k+lM^{k-1}$~for~$l\in[0,K(M)]\cap\bfN.$ By Markov property we see that if~$j_1<j_2\leq M^{k+1},$ then
\begin{eqnarray}\no
H_{j_1} \leq \sup_{x\in [(a-b)j_1^{1/3},~aj_1^{1/3}]}\mbfP^{j_1}_{\L}\left(\forall_{i\leq j_2-j_1},
x+\frac{T_i}{\vartheta}\in[(a-b)(i+j_1)^{\frac{1}{3}},a(i+j_1)^{\frac{1}{3}}]\right)\times H_{j_2}\no
\leq H_{j_2}.
\end{eqnarray}
That is to say, $\{H_{j}\}_{j\leq M^{k+1}}$ is an increasing random sequence of $j.$
We have mentioned above that the choice of $b$ satisfying~$b\in(0,a).$ Combining with~\eqref{e4.3.17}, ~we have
\begin{eqnarray}\label{bu1}
&&\varlimsup\limits_{k\rightarrow +\infty}\frac{\log\left(\sum\limits_{j=M^k+1}^{M^{k+1}}\sup_{\substack{|v'|=j \\V(v')\in\mbr}}\mbfE_{\L}\Big[Z_k^{v'}(\Theta)\Big]\right)}{d^{1/3}_k}\no
\\&\leq&\varlimsup\limits_{k\rightarrow +\infty}\frac{\log\left(\sum\limits_{l=0}^{K(M)}\sum\limits_{j=c_l+1}^{c_{l+1}}e^{\vartheta aM^{\frac{k+1}{3}}-\vartheta(a-b)j^{\frac{1}{3}}}\times H_j\right)}{d^{1/3}_k}\no
\\&\leq&\varlimsup\limits_{k\rightarrow +\infty}\frac{\log\left(\sum\limits_{l=0}^{K(M)}M^{k-1}e^{\vartheta aM^{\frac{k+1}{3}}-\vartheta(a-b)c_l^{\frac{1}{3}}}H_{c_{l+1}}\right)}{d^{1/3}_k}\no
\\&=&\max_{l\in[0, K(M)]}\left[\varlimsup\limits_{k\rightarrow +\infty}\frac{\vartheta aM^{\frac{k+1}{3}}-\vartheta(a-b)c_l^{\frac{1}{3}}}{d^{1/3}_k}+\varlimsup\limits_{k\rightarrow +\infty}\frac{\log H_{c_{l+1}}}{d^{1/3}_k}\right].
\end{eqnarray}
Note that when $M$ is fixed, $K(M)$ is finite and does not depend on $k$, which means that the last equality in \eqref{bu1}. Denote $g_{_M}(x):=\Big(x+\frac{1}{M-1}\Big)^{1/3}.$ By the notation~$c_l:=M^k+lM^{k-1}$ we have
$\varlimsup\limits_{k\rightarrow +\infty}\frac{\vartheta aM^{\frac{k+1}{3}}-\vartheta(a-b)c_l^{\frac{1}{3}}}{d^{1/3}_k}=
\vartheta a g_{_M}(1)-\vartheta(a-b)g_{_M}\Big(\frac{l}{M^2-M}\Big).$
By Corollary \ref{mogc2}, we see
$$\varlimsup\limits_{k\rightarrow +\infty}\frac{\log H_{c_{l+1}}}{d^{1/3}_k}=-\frac{3\gamma_{\sigma}}{\vartheta^2b^2}\left(g_{_M}(1)-g_{_M}\left(\frac{l+1}{M^2-M}\right)\right),~~{\rm \mathbf{P}-a.s.}$$
By the concavity of~$g_{_M}(x),$ we know for any~$l\in[0,K(M)],$ $g_{_M}\left(\frac{l+1}{M^2-M}\right)-g_{_M}\left(\frac{l}{M^2-M}\right)\leq \left(\frac{1}{M^2-M}\right)^{1/3}=\frac{g_{_M}(0)}{M^{1/3}}.$
Hence it is true that
\begin{eqnarray}\label{add1}
&&\varlimsup\limits_{k\rightarrow +\infty}\frac{\log\left(\sum\limits_{j=M^k+1}^{M^{k+1}}\sup_{\substack{|v'|=j \\V(v')\in\mbr}}\mbfE_{\L}\Big[Z_k^{v'}(\Theta)\Big]\right)}{d^{1/3}_k}\no
\\&\leq& \max_{l\in[0, K(M)]}\left[\left(\vartheta a-\frac{3\gamma_{\sigma}}{\vartheta^2b^2}\right)g_{_M}(1) -\left[\vartheta(a-b)-\frac{3\gamma_{\sigma}}{\vartheta^2b^2}\right]g_{_M}\Big(\frac{l}{M^2-M}\Big)+\frac{g_{_M}(0)}{M^{1/3}}\right]\no
\\&\leq& \sup_{x\in[0,1]}\left[\left(\vartheta a-\frac{3\gamma_{\sigma}}{\vartheta^2b^2}\right)g_{_M}(1) +\left(\frac{3\gamma_{\sigma}}{\vartheta^2b^2}-\vartheta(a-b)\right)g_{_M}\left(x\right)+\frac{g_{_M}(0)}{M^{1/3}}\right].
\end{eqnarray}
We observe that the equation~$\frac{3\gamma_{\sigma}}{\vartheta^2b^2}-\vartheta(a-b)=0$ about $b$ has two solutions in $(0,a)$ since~$a>a_c:=\frac{3\sqrt[3]{6\gamma_\sigma}}{2\vartheta}.$ We might as well write them as $b_1$ and $b_2~(b_1<b_2).$
Choose~$b\in(b_1,b_2)$ and note that ~$g_{_M}$ is an increasing function with~$g_{_M}(1)=M^{1/3}g_{_M}(0),$ then we get
\begin{eqnarray}\label{bu1+}
G(M)&:=&\sup_{x\in[0,1]}\left[\left(\vartheta a-\frac{3\gamma_{\sigma}}{\vartheta^2b^2}\right)g_M(1) +\left(\frac{3\gamma_{\sigma}}{\vartheta^2b^2}-\vartheta(a-b)\right)g_M\left(x\right)+\frac{g_M(0)}{M^{1/3}}\right]\no
\\&=&\left[\left(\vartheta a-\frac{3\gamma_{\sigma}}{\vartheta^2b^2}\right)M^{1/3} +\frac{3\gamma_{\sigma}}{\vartheta^2b^2}+\vartheta b-\vartheta a+\frac{1}{M^{1/3}}\right]g_M(0).
\end{eqnarray}

Recall that $r_k=e^{k^{v}},~v\in(\frac{2}{\lambda_3},\frac{1}{3}).$ Hence when $b\in(b_1, b_2),$ we have
\begin{eqnarray}\label{add2}
&&\varlimsup\limits_{k\rightarrow +\infty}\frac{\log\left(1+(r_k-1)\sum\limits_{j=M^k+1}^{M^{k+1}}\sup_{\substack{|v'|=j \\V(v')\in\mbr}}\mbfE_{\L}\Big[Z_k^{v'}(\Theta)\Big]\right)}{d^{1/3}_k}\no
\\&=&\varlimsup\limits_{k\rightarrow +\infty}\frac{\log 1}{d^{1/3}_k}\vee\varlimsup\limits_{k\rightarrow +\infty}\frac{\log\left((r_k-1)\sum\limits_{j=M^k+1}^{M^{k+1}}\sup_{\substack{|v'|=j \\V(v')\in\mbr}}\mbfE_{\L}\Big[Z_k^{v'}(\Theta)\Big]\right)}{d^{1/3}_k}\no
\\&\leq&G(M).
\end{eqnarray}

In the light of $\eqref{impc4},$ it is time to estimate the lower bound of $\mbfE_{\L}(Z_{k,b}).$ Recall~$c_0:=M^k, V(\tilde{z}) =ac_0^{\frac{1}{3}}-\vartheta^{-1} K_{c_0}$ and $d_k:=M^{k+1}-M^{k}.$ 
For any~$\epsilon\in(0,b),$ by the definition of~$Z_{k,b},~ \eqref{shift}$ and \eqref{mto0} we can see
\begin{eqnarray}\label{e4.3.31} \mbfE_{\L}Z_{k,b}&=&\mbfE_{\L}\Bigg(\sum_{|u|=M^{k+1}}\1_{\{u\in\Theta\}}\Bigg)\no
\\&=&\mbfE^{c_0}_{\L}\Bigg(\sum_{|u|=d_k}\1_{\left\{\substack{\forall i\leq d_k, N(u_{i-1})\leq r_{k}, V(u_i)+\vartheta^{-1}K_{c_0+i}+V(\tilde{z})\in[(a-b)(i+c_0)^{\frac{1}{3}},~a(i+c_0)^{\frac{1}{3}}]}\right\}}\Big|V(u_0)=0\Bigg)\no
\\&=&\mbfE^{c_0}_{\L}\Big(e^{T_{d_k}}\1_{\{0<i\leq d_k,\xi_i\leq r_k,T_i\in[\vartheta(a-b)(i+c_0)^{\frac{1}{3}}-\vartheta ac^{\frac{1}{3}}_0,\vartheta a(i+c_0)^{\frac{1}{3}}-\vartheta ac^{\frac{1}{3}}_0]\}}|T_0=0\Big)\no
\\&\geq&e^{\vartheta (a-\epsilon)M^{\frac{k+1}{3}}-\vartheta ac_0^{\frac{1}{3}}}\mbfP^{c_0}_{\L}\left(\substack{\forall 0<i\leq d_k,T_i\in\left[\vartheta(a-b)(i+c_0)^{\frac{1}{3}},\vartheta a(i+c_0)^{\frac{1}{3}}\right]\\~~\\\xi_i\leq r_k,~
T_{d_k}\in\Big[\vartheta(a-\epsilon)M^{\frac{k+1}{3}},\vartheta aM^{\frac{k+1}{3}}\Big]}\Bigg|T_0=ac^{\frac{1}{3}}_0\right).\end{eqnarray}
 Applying the Corollary \ref{mog0} \footnote{Note that we can not utilize~Corollary \ref{mogc2} here for the reason of~$T_0=\vartheta a(0+c_0)^{\frac{1}{3}},$ which is located on the boundary but not the interior of the interval $\left[\vartheta(a-b)i^{1/3},\vartheta ac_0^{1/3}\right]$. That is why we need Condition 4 in Section 2.}, we obtain
 $$ \varliminf\limits_{k\rightarrow+\infty}\frac{\log \mbfE_{\L}(Z_{k,b})}{d^{1/3}_k}\geq \vartheta(a-\epsilon)\Big(\frac{M}{M-1}\Big)^{\frac{1}{3}}-a\vartheta\Big(\frac{1}{M-1}\Big)^{\frac{1}{3}}-\frac{\gamma_\sigma}{b^2\vartheta^2}\int_{0}^1\Big(x+\frac{1}{M-1}\Big)^{-\frac{2}{3}}dx. $$
Letting $\epsilon\downarrow0$, we get
 \begin{eqnarray}\label{e4.3.32} \varliminf\limits_{k\rightarrow+\infty}\frac{\log \mbfE_{\L}(Z_{k,b})}{d^{1/3}_k}&\geq&
\Big(a\vartheta-\frac{3\gamma_\sigma}{b^2\vartheta^2}\Big)\big(g_M(1)-g_M(0)\big)\no
\\&=&\Big(a\vartheta-\frac{3\gamma_\sigma}{b^2\vartheta^2}\Big)(M^{1/3}-1)g_M(0),~~{\rm \mathbf{P}-a.s.}\end{eqnarray}
From the above discussion and the fact $a\vartheta-\frac{3\gamma_\sigma}{b^2\vartheta^2}>0$ we can see $\lim\limits_{k\rightarrow+\infty}\mbfE_{\L}(Z_{k,b})=+\infty,~{\rm \mathbf{P}-a.s.},$ which means that $\lfloor\eta\mbfE_{\L}(Z_{k,b})\rfloor\geq \frac{\eta}{2}\mbfE_{\L}(Z_{k,b})$ for large enough $k.$ Then~\eqref{e4.3.32} tells us
\begin{eqnarray*}\varliminf\limits_{k\rightarrow+\infty}\frac{\log \lfloor\eta\mbfE_{\L}(Z_{{k-1},b})\rfloor}{d^{1/3}_k}&\geq&
\frac{1}{M^{1/3}}\Big(a\vartheta-\frac{3\gamma_\sigma}{b^2\vartheta^2}\Big)(M^{1/3}-1)g_M(0),~~{\rm \mathbf{P}-a.s.}\end{eqnarray*}
Recall the definition of $A_{k,\L}$ in \eqref{AKL}. Combining $\eqref{e4.3.32}$ with~$\eqref{add2}$ we get
\begin{eqnarray}\label{e4.3.33} &&\varliminf\limits_{k\rightarrow+\infty}\frac{\log\left(\lfloor\eta\mbfE_{\L}(Z_{k-1,b})\rfloor A_{k,\L}\right)}{d^{1/3}_k}\no
\\&\geq&\Big(a\vartheta-\frac{3\gamma_\sigma}{b^2\vartheta^2}\Big)\left(1-\frac{1}{M^{1/3}}\right)g_M(0)
-G(M)\no
\\&=&\left[\left(\vartheta a-\vartheta b-\frac{3\gamma_{\sigma}}{\vartheta^2b^2}\right)-\frac{2}{M^{1/3}}\right]g_M(0).
\end{eqnarray}
Note that $\vartheta a>\frac{3\gamma_\sigma}{b^2\vartheta^2}+\vartheta b$ for $b\in(b_1,b_2).$ We choose a large enough constant $M$ such that
$\left(\vartheta a-\vartheta b-\frac{3\gamma_{\sigma}}{\vartheta^2b^2}\right)-\frac{2}{M^{1/3}}>0,$ which means
$$\varliminf\limits_{k\rightarrow+\infty}\frac{\log\left(\lfloor\eta\mbfE_{\L}(Z_{k-1,b})\rfloor A_{k,\L}\right)}{d_k}>0,~{\rm \mathbf{P}-a.s.}$$
Thus $\eqref{e4.3.5}$ holds. Recalling the analysis at the beginning of this proof we finally get~\eqref{e4.3.5+}. So far we have shown $\mbfP_{\L}(\mathcal{S})>0, ~{\rm \mathbf{P}-a.s.}$ when the barrier function with parameter $\alpha=\frac{1}{3}$ and $a>a_c.$ \qed

\noindent{\bf Proof of Theorem \ref{2} (2d)}

By reviewing the above proof again, we can even get Theorem \ref{2} (2d). According to $\eqref{e4.3.5--},$ we get~$\mbfP_{\L}(\mathcal{S})>0, ~{\rm \mathbf{P}-a.s.}$ by proving
~$\lim\limits_{n\rightarrow +\infty}P_{n,\L}>0,~{\rm \mathbf{P}-a.s.}$ By the definition of~$P_{n,\L}$ we see
\begin{eqnarray}\label{e4.3.35}
\lim\limits_{n\rightarrow +\infty}P_{n,\L}&=&\mbfP_{\L}\left(\forall k\in\bfN,\sharp\{u\in\mathcal{T}_{M^{k}}, \forall i\leq M^{k}, V(u_i)\leq ai^{\frac{1}{3}}-\vartheta^{-1}K_{i}\}\geq\eta\mbfE_{\L}(Z_{k-1,b})\right)\no
\\&=&\mbfP_{\L}\left(\forall k\in\bfN,Y_{M^k}\geq\eta\mbfE_{\L}(Z_{k-1,b})\right).
\end{eqnarray}
According to $\eqref{e4.3.32}$ and the fact $a\vartheta-\frac{3\gamma_\sigma}{b^2\vartheta^2}>b\vartheta,$ we obtain
$\varliminf\limits_{k\rightarrow+\infty}\frac{\log \mbfE_{\L}(Z_{k,b})}{d^{1/3}_k}\geq
b\vartheta\big(g_M(1)-g_M(0)\big).$ Hence we have $g_M(1)\rightarrow 1, ~g_M(0)\rightarrow 0,\frac{d_{k-1}}{M^{k}}\rightarrow 1$ as~$M\rightarrow +\infty.$
 It means that for any $~\varepsilon>0,$ we can find a large enough~$M$ such that $\varliminf\limits_{k\rightarrow+\infty}\frac{\log (\eta\mbfE_{\L}(Z_{{k-1},b}))}{M^{k/3}}\geq
b\vartheta-\varepsilon,~{\rm \mathbf{P}-a.s.}$ Note that~$b\in(b_1,b_2),$ then by~\eqref{e4.3.35} we have
$$ \mbfP_{\L}\left(\varliminf\limits_{k\rightarrow+\infty}\frac{\log Y_{M^k}}{M^{k/3}}\geq b_2\vartheta-\varepsilon\right)\geq\lim\limits_{n\rightarrow +\infty}P_{n,\L}>0,~~{\rm \mathbf{P}-a.s.},$$
which is the conclusion in~$\eqref{big1}.$\qed 


At last, we turn to~Theorem \ref{1} (1a).

\noindent{\bf The proof of~Theorem \ref{1} (1a)}

Let~$a_{c+}$ be a constant such that ~$a_{c+}>a_c$. Define~$j_n:=(a_{c+})^{\frac{1}{\alpha}}n^{\frac{1}{3\alpha}}-a^{\frac{1}{\alpha}}n,~n\in\bfN^+.$ Since the case in Theorem \ref{1} (1a) is~$\alpha>1/3,$ we have~$j_{max}:=\max_{n\in\bfN^+}j_n<+\infty.$ Choose $k$ large enough such that~$ka^{\frac{1}{\alpha}}>j_{max},$
 which ensures that~$a(n+k)^{\alpha}>a_{c+}n^{1/3}, \forall n\in\bfN^+.$ Note that $\alpha>1/3,$ hence it is true that~$\inf_{n\in\bfN^+}(a(n+k)^{\alpha}-a_{c+}n^{1/3})>0.$ we can find $a_->0$ small enough such that
~$a(n+k)^{\alpha}>a_-k+a_{c+}n^{1/3}, \forall n\in\bfN^+$ and~$a_-<\min\{ak^{\alpha-1},a\}.$ In this way we can ensure that
~$ai^{\alpha}> a_-i$ for $1\leq i\leq k$
~and~$ai^{\alpha}>a_-k+a_{c+}(i-k)^{1/3}$ for $i>k.$ By Markov property we see
\begin{eqnarray*}
\mbfP_{\L}(\mathcal{S})&=&\mbfP_{\L}(\exists u\in \mathcal{T_{\infty}}, \forall i\in\bfN,  V(u_i)\leq ai^{\alpha}-\vartheta^{-1}K_i)
\\&\geq& \mbfP_{\L}(\exists u\in \mathcal{T}_k, \forall i\leq k,  V(u_i)\leq a_-i-\vartheta^{-1}K_i)
\\&\times&\mbfP^k_{\L}(\exists u\in \mathcal{T}^z_{\infty}, \forall i\in \bfN,  V(u_i)\leq a_{c+}i^{1/3}+\vartheta^{-1}(K_{k+i}-K_{k})|V(z)=0)
\\&:=&U_1\times U_2,
\end{eqnarray*}
where~$\mathcal{T}^z_{\infty}$ represents a infinite path in~$\mathcal{T}^z$ and $\mathcal{T}^z$ the genealogical tree with ancestor~$z.$ Theorem \ref{1} (1b) tells us~$U_2>0, {\rm \mathbf{P}-a.s.}$ and Theorem \ref{2} (2c)~means that~$U_1>0, ~{\rm \mathbf{P}-a.s.}$ Hence we have~$\mbfP_{\L}(\mathcal{S})>0,~{\rm \mathbf{P}-a.s.}$ \qed


 \ack
We would like to thank the referees greatly for their careful review and valuable suggestions.
This work is supported by the Fundamental Research Funds for the Central Universities (NO.2232021D-30) and the National Natural Science Foundation of China (NO.11971062).



\begin{thebibliography}{99}
\bibitem{A2013}
A\"{\i}d\'{e}kon, E.~(2013)~Convergence in law of the minimum of a branching random walk. {\it Ann. Probab.} {\bf 41(3A)}, 1362-1426.
\bibitem{AJ2011}
A\"{\i}d\'{e}kon, E. and Jaffuel, B.~(2011)~Survival of branching random walks with absorption. {\it Stochastic Proc. Appl.} {\bf 121}, 1901-1937.
\bibitem{AR2009}
Addario-Berry, L. and Reed, B.~(2009)~Minima in branching random walks. {\it Ann. Probab.} {\bf 37}, 1044-1079.
\bibitem{BC1993}
Baillon, J.B., Cl\'{e}ment, P., Greven, A. and Hollander, F.~(1993)~A variational approach to branching random walk in
random environment. {\it Ann. Probab.} {\bf 21}(1), 270-317.
\bibitem{BG2011}
B\'{e}rard, J., Gou\'{e}r\'{e},J.B.~(2011)~Survival probability of the branching random
walk killed below a linear boundary. {\it Electronic. J. Probab.} 16(14): 396-418.
\bibitem{B1976}
Biggins, J. D.~(1976)~The first and last birth problems for a multitype age-dependent branching process. {\it Adv. Appl. Probab.} {\bf 8}, 446-459.
\bibitem{BK2004}
Biggins, J.D. and Kyprianou, A.E.~(2004)~Measure change in multitype branching. {\it Adv. Appl. Probab.} {\bf 36}(2), 544-581.
\bibitem{BLSW1991}
Biggins, J.D., Lubachevsky,B.D., Shwartz,A. and Weiss,A.~(1991)~A branching random walk with a barrier. {\it Ann. Appl.
Probab.} {\bf 1}, 573-581.
\bibitem{CR1979}
Cs\"{o}rg\H{o}, M. and R\'{e}v\'{e}sz, P.~(1979)~How big are the increments of a Wiener process?~
{\it Acta Mathematica Academiae Scientiarum Hungarica.} {\bf 33}(1-2), 37-49.
\bibitem{DS2007}
Derrida, B. and Simon, D.~(2007)~The survival probability of a branching random walk in presence of an absorbing wall. {\it Europhys. Lett.} {\bf 78}(6), 346-350.
\bibitem{DS2008}
Derrida, B. and Simon, D.~(2008)~Quasi-stationary regime of a branching random walk in presence of an absorbing wall. {\it J. Stat. Phys.} {\bf 131}, 203-233.
\bibitem{GHS2011}
Gantert, N., Hu, Y. and Shi, Z.~(2011)~Asymptotics for the survival probability in a killed branching random walk. {\it Ann. Inst. Henri Poincar\'{e} Probab. Stat.} {\bf 47}(1), 111-129.
\bibitem{GLW2014}
Gao, Z. Liu, Q. and Wang, H.~(2014)~Central limit theorems for a branching random walk with a random environment in time.
{\it Acta Math. Sci. Ser. B Engl. Ed.} {\bf 34}(2), 501-512.
\bibitem{GL2016}
Gao, Z. and Liu, Q.~(2016)~Exact convergence rates in central limit theorems for a branching random walk with a random environment in time.
{\it Stochastic Proc. Appl.} {\bf 126}(9), 2634-2664.
\bibitem{H1974}
Hammersley, J. M.~(1974)~Postulates for subadditive processes. {\it Ann. Probab.} {\bf 2}, 652-680.
\bibitem{HH2007}
Harris, J. W. and Harris, S. C.~(2007)~Survival probabilities for branching Brownian motion with absorption. {\it Electron. Commun. Probab.} {\bf 12}, 81-92.
\bibitem{HL2014}
Huang, C. and Liu, Q.~Branching random walk with a random environment in time. ArXiv:1407.7623.
\bibitem{HY2009}
Hu, Y. and Yoshida, N.~(2009)~Localization for branching random walks in random environment. {\it Stochastic Proc. Appl.} {\bf 119}(5), 1632-1651.
\bibitem{HS2009}
Hu, Y. and Shi, Z.~(2009)~Minimal position and critical martingale convergence in branching random walks, and
directed polymers on disordered trees. {\it Ann. Probab.} {\bf 37}(2), 742-789.
\bibitem{K1978}
Kesten, H.~(1978)~Branching Brownian motion with absorption. {\it Stochastic Proc. Appl.} {\bf 7}, 9-47.
\bibitem{BJ2012}
Jaffuel, B.~(2012)~The critical barrier for the survival of branching random walk with absorption. {\it Ann. Inst. Henri Poincar\'{e} Probab. Stat.} {\bf 48}(4), 989-1009.
\bibitem{KP1976}
Kahane, J. P. and Peyri\'{e}re, J.~(1976)~Sur certaines martingales de Benoit Mandelbrot. {\it Adv. Math.} {\bf 22}(2), 131-145.
\bibitem{K1975}
Kingman, J. F. C.~(1975)~The first birth problem for an age dependent branching process. {\it Ann. Probab.} {\bf 3}, 790-801.
\bibitem{LSW1989}
Lubachevsky, B., Shwartz, A. and Weiss, A.~(1989)~Rollback sometimes works ... if
filtered. {\it In Proceedings of the Winter Simulation Conference} IEEE, New York. 630-639.
\bibitem{LSW1990}
Lubachevsky, B., Shwartz, A. and Weiss, A.~(1990)~\emph{An analysis of rollback-based simulation.} EE PUB 755, Technion, Israel.
\bibitem{LZ2019}
Liu, J. and Zhang, M.~(2019)~Critical survival barrier for branching random walk. {\it Front. Math. China} {\bf 14}(6), 1259-1280.
\bibitem{LY201801}
Lv, Y.~(2019)~Brownian motion between two random trajectories. {\it Markov Process. Related Fields} {\bf 25}(2), 359-377.
\bibitem{Lv201802}
Lv, Y. and Hong, W.~(2023)~Quenched small deviation for the trajectory of a random walk with random environment in time. {\it Theory Probab. Appl.} {\bf 68}(2), 322-343.
\bibitem{M2015a}
Mallein,~B.~(2015)~Maximal displacement in a branching random walk through interfaces. {\it Electron. J. Probab.} {\bf 68}(20), 1-40.
\bibitem{M2015b}
Mallein,~B.~(2015)~Maximal displacement of a branching random walk in time-inhomogeneous environment. {\it Stochastic Process. Appl.} {\bf 125}, 3958-4019.
\bibitem{M2017}
Mallein, B.~(2017)~Branching random walk with selection at critical rate. {\it Bernoulli.} {\bf23}(3), 1784-1821.
\bibitem{MM2016}
Mallein, B. and Mi{\l}o\'{s}, P.~(2019)~Maximal displacement of a supercritical branching random walk in a time-inhomogeneous random environment. {\it Stochastic Process. Appl.} {\bf129}, 3239-3260.
\bibitem{Mog1974}
Mogul'ski\v{\i}, A. A.~(1974)~Small deviations in the space of trajectories. {\it Theory Probab. Appl.} {\bf 19}, 726-736.
\bibitem{P1974}
Peyri\'{e}re, J.~(1974)~Turbulence et dimension de Hausdorff. {\it C. R. Acad. Sci. Paris S\'{e}r. A}. {\bf 278}, 567-569.
\bibitem{Sak2006}
Sakhanenko, A.I. (2006) Estimates in the invariance principle in terms of truncated power moments. {\it Siberian Math. J.} {\bf47}(6), 1113-1127.
\bibitem{S2015}
Shi, Z.~(2015)~{\it Branching random walks.} \'{E}cole d'\'{E}t\'{e} de Probabilit\'{e}s de Saint-Flour XLII-2012. Lecture Notes in Mathematics 2151, ~Springer, Berlin.
\bibitem{WH2017}
Wang, X. and Huang, C.~(2017)~Convergence of Martingale and Moderate Deviations for a Branching Random Walk with a Random Environment in Time.
{\it J. Theoretical Probab.} {\bf 30}, 961-995.
\bibitem{Z2004}
Zeitouni, O.~(2008)~Part II: \emph{Random Walks in Random Environment.} Lectures on Probability Theory and Statistics. Springer, Berlin, Heidelberg.
\end{thebibliography}
\end{document}